%% file: main.tex
\documentclass[10pt,dvipsnames]{article}
\usepackage{geometry}
\usepackage{mathtools}
\usepackage{graphicx}
\usepackage[pdftex,colorlinks=true]{hyperref}
\usepackage{caption}
\usepackage{physics}
\usepackage{tensor}
\usepackage[most]{tcolorbox}
\tcbuselibrary{theorems}
\usepackage{empheq}
\usepackage{enumitem}
\usepackage{wrapfig}
\usepackage{pythonhighlight}

\usepackage{float}
\usepackage{graphicx}

\input{customstyle.tex}

\makeatletter
\newcommand{\linebreakand}{%
  \end{@authorhalign}
  \hfill\mbox{}\par
  \mbox{}\hfill\begin{@authorhalign}
}
\makeatother

\title{High-dimensional robust regression under heavy-tailed data: Asymptotics and Universality}
\author{Urte Adomaityte\\\small{Department of Mathematics}\\\small{King's College London}\\
\small{United Kingdom} 
\and 
Leonardo Defilippis\\\small{Departement d'Informatique}\\
\small{\'Ecole\,Normale\,Sup\'erieure}
\\\small{PSL \& CNRS}\\
\small{France}
\and 
Bruno Loureiro\\\small{Departement d'Informatique}\\ 
\small{\'Ecole\,Normale\,Sup\'erieure}
\\\small{PSL \& CNRS}\\
\small{France} 
\and
Gabriele Sicuro\\\small{Department of Mathematics}\\\small{University of Bologna}\\\small{Italy}}

\date{\today}

\allowdisplaybreaks
\begin{document}

\maketitle
\begin{abstract}
We investigate the high-dimensional properties of robust regression estimators in the presence of heavy-tailed contamination of both the covariates and response functions. In particular, we provide a sharp asymptotic characterisation of M-estimators trained on a family of elliptical covariate and noise data distributions including cases where second and higher moments do not exist. We show that, despite being consistent, the Huber loss with optimally tuned location parameter $\delta$ is suboptimal in the high-dimensional regime in the presence of heavy-tailed noise, highlighting the necessity of further regularisation to achieve optimal performance. This result also uncovers the existence of a transition in $\delta$ as a function of the sample complexity and contamination. Moreover, we derive the decay rates for the excess risk of ridge regression. We show that, while it is both optimal and universal for covariate distributions with finite second moment, its decay rate can be considerably faster when the covariates' second moment does not exist. Finally, we show that our formulas readily generalise to a richer family of models and data distributions, such as generalised linear estimation with arbitrary convex regularisation trained on mixture models.
\end{abstract}

\section{Introduction}
\label{sec:intro}
\input{sections/introduction}
\section{High-dimensional asymptotics} 
\input{sections/asymptotics}
\section{Discussion}
\label{sec:examples}
\input{sections/discussion}

\section{Generalised linear estimation on an elliptical mixture}
\label{sec:gencluster}
\input{sections/generalisation}
\subsubsection*{Acknowledgements} 
\input{sections/acknowledgement}

\appendix
\section*{\LARGE Appendix}
\section{Replica derivation of the fixed-point equations}
\label{app:replica}
\input{sections/appendix/replicas}
\section{Asymptotic results for ridge-regularised losses}
\label{app:ridge_regul}
\input{sections/appendix/ridge}

\section{Further numerical results}
\label{app:furthernumerics}
\input{sections/appendix/numerics}



\bibliographystyle{plainnat}
\bibliography{refs.bib}
\end{document}

%% file: customstyle.tex
\usepackage{amsmath,amssymb,mathrsfs,bm,mathtools,amsthm}
\usepackage{latexsym,amscd,amsbsy,dsfont,amsfonts,xfrac,upgreek}
\usepackage{physics,tensor,empheq,enumitem}
\usepackage{xcolor}
\usepackage{tikz}
\usepackage{adjustbox}
\setlist{nosep}
\usepackage{mlmodern}
\usepackage[T1]{fontenc}
\usepackage{hyperref}
\hypersetup{pdfauthor={ALDS},pdftitle={Robust},%
            colorlinks, linktocpage=true, pdfstartpage=1, pdfstartview=FitV,%
    breaklinks=true, pdfpagemode=UseNone, pageanchor=true, pdfpagemode=UseOutlines,%
    plainpages=false, bookmarksnumbered, bookmarksopen=true, bookmarksopenlevel=1,%
    hypertexnames=true, pdfhighlight=/O,%
    urlcolor=orange, linkcolor=blue, citecolor=ForestGreen, 
        }

\newcommand{\comprimi}{\medmuskip=0mu
\thinmuskip=0mu
\thickmuskip=0mu}

\DeclareMathOperator{\e}{e}
\DeclareMathOperator{\argmin}{argmin}

\usepackage[numbers]{natbib}

\DeclareMathOperator{\erfc}{erfc}

\newcommand{\R}{{\mathds{R}}}

\DeclareMathOperator*{\Extr}{Extr}

\newcommand{\sq}{{{\mathsf q}}}

\newcommand{\shq}{{\hat{\mathsf q}}}

\newcommand{\be}{{\boldsymbol{e}}}
\newcommand{\bx}{{\boldsymbol{{x}}}}
\newcommand{\bz}{{\boldsymbol{{z}}}}

\newcommand{\bmu}{{\boldsymbol{{\mu}}}}

\newcommand{\bbeta}{{\boldsymbol{{\beta}}}}
\newcommand{\bEta}{{\boldsymbol{{\eta}}}}

\newcommand{\bxi}{{\boldsymbol{{\xi}}}}

\newcommand{\bI}{{\boldsymbol{{I}}}}
\newcommand{\bUno}{{\boldsymbol{{1}}}}
\newcommand{\bM}{{\boldsymbol{{M}}}}

\newcommand{\bt}{{\boldsymbol{{t}}}}
\newcommand{\bht}{{\boldsymbol{\hat{t}}}}
\newcommand{\bQ}{{\boldsymbol{{Q}}}}

\newcommand{\bgg}{{\boldsymbol{{g}}}}

\newcommand{\hQ}{{\hat Q}}
\newcommand{\hV}{{\hat V}}
\newcommand{\hR}{{\hat R}}

\newcommand{\bhM}{{\boldsymbol{{\hat M}}}}

\newcommand{\bhQ}{{\boldsymbol{{\hat Q}}}}

\newcommand{\bSigma}{{\boldsymbol{{\Sigma}}}}

\newcommand{\hsigma}{{\hat\sigma}}
\newcommand{\tsigma}{{\tilde\sigma}}

\makeatletter
\newsavebox{\@brx}
\newcommand{\llangle}[1][]{\savebox{\@brx}{\(\m@th{#1\langle}\)}%
  \mathopen{\copy\@brx\kern-0.5\wd\@brx\usebox{\@brx}}}
\newcommand{\rrangle}[1][]{\savebox{\@brx}{\(\m@th{#1\rangle}\)}%
  \mathclose{\copy\@brx\kern-0.5\wd\@brx\usebox{\@brx}}}
\makeatother

\newtheorem{thm}{Theorem}[section]
\newtheorem{result}[thm]{Result}
\newtheorem{assumption}[thm]{Assumption}

%% file: sections/introduction.tex
Consider the classical statistical problem of estimating a vector $\bbeta_{\star}\in\R^{d}$ from $n$ i.i.d.~pairs of observations $\mathcal{D}\coloneqq\{(\bx_{i}, y_{i})\in\R^{d+1}: i\in [n]\}$ from a ``teacher'' probabilistic model
\begin{subequations}
\label{eq:def}    
\begin{equation}
    y_{i} \sim P_0(\cdot|\bbeta_{\star}^\intercal\bx_{i}),\label{eq:labels}
\end{equation}
where $\bx_{i}\in\R^{d}$ are the covariates. A relevant special case is the linear setting,
\begin{equation}
\label{eq:def:data}
    y_{i} =\bbeta_{\star}^\intercal\bx_{i}+\eta_i,
\end{equation}
where $\eta_i\in\R$ is some label noise, which we assume to be a random quantity with zero mean and independent of the covariates This manuscript is concerned with the characterisation of the following class of (regularised) M-estimators:
\begin{equation}
\label{eq:def:erm}
    \hat{\bbeta}_{\lambda}\coloneqq \underset{\bbeta\in\R^{d}}{\argmin}\sum\limits_{i=1}^{n}\rho(y_{i}-\bbeta^\intercal \bx_{i})+\frac{\lambda}{2}\|\bbeta\|_2^2,\qquad \lambda\in\R_+\coloneqq [0,+\infty),
\end{equation}
\end{subequations}
with $\rho\colon\R\to\R_{+}$ a convex objective function. A popular example is least-squares regression, where $\rho(t)=\sfrac{t^2}{2}$ and $\lambda=0$; in this case, assuming a linear teacher as in Eq.~\eqref{eq:def:data}, $\hat{\bbeta}_{\lambda}$ is the maximum likelihood estimator for $\bbeta_\star$ if $\eta_i\sim\mathcal{N}(0,1)$. It is well-known, however, that the least-squares estimator suffers in the presence of outliers in the data \citep{Huber1973}. Indeed, the fact that the gradient of the loss $\rho'(t) = t$ is unbounded implies that an outlier can have arbitrary influence over the solution of Eq.~\eqref{eq:def:erm}. Tailoring the objective $\rho$ to be insensitive (i.e., \emph{robust}) to outliers in the training data $\mathcal{D}$ is a classical statistical problem \citep{huber2004robust, hampel2011robust, rousseeuw2005robust, maronna2019robust}. In his seminal work, \citet{Huber64} has shown that judiciously trimming the squared loss
\begin{align}
\label{eq:def:huberloss}
    \rho_{\delta}(t) = 
    \begin{cases}
        \sfrac{t^2}{2} & \text{ if } |t|<\delta\\
        \delta |t|-\sfrac{\delta^2}{2} & \text{otherwise}
    \end{cases}
\end{align}
for $\delta\geq 0$ provides a convenient solution to this problem while preserving the convexity of the task in Eq.~\eqref{eq:def:erm}. Indeed, besides enjoying standard statistical guarantees in the classical limit of $n\to\infty$, such as consistency and asymptotic normality \citep{Huber1973, van2000asymptotic}, the so-called \emph{Huber loss} in Eq.~\eqref{eq:def:huberloss} has been shown to be optimal in different regards. For instance, \citet{Huber64} has shown it achieves minmax asymptotic variance under symmetric contamination of the normal distribution. It also has the smallest asymptotic variance among losses with bounded sensitivity to outliers, as it can formalised by Hampel's influence function \citep{Hampel1974}. However, while these guarantees are fit for a classical statistical regime where data is abundant $(n\gg d)$, they fall short in modern tasks where the number of features can be comparable to the quantity of data ($n\approx d$). Investigating the properties of estimators in the proportional \emph{high-dimensional regime} where $n,d\to\infty$ at fixed sample complexity $\alpha=\sfrac{n}{d}$ has been a major endeavor in the statistical literature in the past decade, where it has been shown that standard guarantees for the maximum likelihood estimator, such as unbiasedness \citep{Javanmard2018, SurCandes2019, Sur2020, Bellec2022b, ZhaoSurCandes2022} and calibration \citep{Bai21,Clarte2023a,Clarte2023b} break down in this regime. However, the majority of these works have focused on the case of sub-Gaussian features and bounded noise variance. Our goal in this manuscript is to go beyond this assumption by providing a high-dimensional characterisation of M-estimators for a family of heavy-tailed distributions for both the covariates and noise distributions. 

\paragraph{Heavy-tailed data ---}  In the following, we consider a family of covariate distributions parametrised as
\begin{equation}
\label{eq:def:superstats}
    \bx_i=\sigma_i\bz_i,\qquad \bz_i\sim\mathcal N(\mathbf 0,\sfrac{1}{d}\bI_d)\qquad \sigma_i\sim\varrho,
\end{equation}
where for each $i\in[n]$, $\sigma_i\in \R^*_+\coloneqq (0,+\infty)$ is an independent random variable with probability density $\varrho$ supported on the positive real line. This class of covariates appeared under different contexts in non-equilibrium statistical physics \citep{Beck2003}, statistics \citep{Karoui2018b,adomaityte2023} natural image modelling \citep{Wainwright1999}, hierarchical data priors in Bayesian modelling \cite{Gelman_Hill_2006}, quantitative finance \cite{Delpini_Bormetti15}, where it has been shown that by suitably choosing $\varrho$ yields a large family of \textit{power-law tailed distributions} for $\bx$, see Table \ref{tab:examples} for concrete examples. In particular, we are interested in investigating the impact of heavy-tail contamination of both the covariates and responses, and the sensitivity of M-estimators of the type in Eq.~\eqref{eq:def:erm} to them. More concretely, in the following we consider the \emph{Huber $\epsilon$-contamination model} for the covariates by assuming
\begin{equation}
\label{eq:def:hubercont}
    p(\bx)=\mathbb E_\sigma[\mathcal N(\bx;\mathbf 0,\sfrac{\sigma^2}{d}\bI_d)]\quad\text{with}\quad \sigma\sim\varrho=(1-\epsilon_{\rm c})\delta_{\sigma,1}+\epsilon_{\rm c}\varrho_0,\quad\epsilon_{\rm c}\in[0,1],
\end{equation} 
where $\varrho_0$ is a density over $\R^*_+$, and $\delta_{\sigma,1}$ equals to $1$ when $\sigma=1$ and is zero everywhere else. The distribution $p(\bx)$ belongs therefore to the family defined by Eq.~\eqref{eq:def:superstats}: the quantity $\epsilon_{\rm c}$ measures, therefore, the \textit{contamination} of a pure Gaussian covariate distribution (recovered for $\epsilon_{\rm c}=0$) with a possibly heavy-tailed law. Moreover, our results can be extended to the case of \textit{mixtures}, i.e., the case of the covariates grouped in $K$ different clouds. This generalisation is introduced and briefly discussed in Section~\ref{sec:gencluster}.

\begin{table}
\centering
\begin{adjustbox}{max width=\textwidth}
\def\arraystretch{1.6}
\begin{tabular}{| l c | c | c |}
\hline
&$\varrho(\sigma)$& $p(x) = \mathbb E_\sigma[\mathcal N(x;0,\sigma^2)]$ & $k$th moment exists if \\ \hline\hline 
Inverse-Gamma\hspace{1cm}  & $\frac{2b^{a}\exp(-\sfrac{b}{\sigma^2})}{\Gamma(a)\sigma^{2a+1}}$ & $\frac{(2b)^a \Gamma(a+\sfrac{1}{2})}{\sqrt{\pi}\Gamma(a) (2b+x^2)^{a+\sfrac{1}{2}}}$ & $k< 2a$ \\ 
  \multicolumn{2}{|r|}{for $a=b=\sfrac{1}{2}$} 
  & $\frac{1}{\pi}\frac{1}{1+x^2}$ (Cauchy) & none \\ 
  \multicolumn{2}{|r|}{for $2a=2b=n$} 
  & $\comprimi\frac{\Gamma(\frac{n+1}{2})}{\sqrt{\pi n}\Gamma(\sfrac{n}{2})}\big(1+\frac{x^2}{n}\big)^{-\frac{n+1}{2}}$ (Student-t)& $k<n$ \\

  \multicolumn{2}{|r|}{for $a=1+b\to+\infty$} 
  & $\frac{1}{\sqrt{2\pi}}\e^{-\frac{x^2}{2}}$ (Gaussian)& all \\
\hline
Pareto& $\frac{2a \theta(\sigma-1)}{\sigma^{2a+1}}$ & $\frac{a 2^a\gamma(a+\sfrac{1}{2},\sfrac{1}{2}x^2)}{\sqrt\pi |x|^{2a+1}} $  & $k<2a$  \\
  \multicolumn{2}{|r|}{for $a\to+\infty$}
  & $\frac{1}{\sqrt{2\pi}}\e^{-\frac{x^2}{2}}$ (Gaussian) & all \\
\hline
\end{tabular}
\end{adjustbox}
\caption{Concrete examples of distributions $p(x)=\mathbb E_\sigma[\mathcal N(x;0,\sigma^2)]$ for different densities $\varrho$ of $\sigma$. Here $\gamma(a,x)$ is the lower incomplete Gamma function. Note that relevant distributions, such as the Cauchy and the Student's $t$-distribution, appear as special cases. The last column shows the values of the parameters of $\varrho$ for which the $k$th moment of the distribution is finite.}
\label{tab:examples}
\end{table}

\paragraph{Main contributions ---}  The \textbf{key contributions} in this manuscript are:
\begin{itemize}[noitemsep,topsep=0pt]
\item We provide an asymptotic characterisation of the statistics of the M-estimator $\hat\bbeta_\lambda$ defined by Eqs.~\eqref{eq:def} with heavy-tailed covariates and general label noise distributions in the proportional high-dimensional regime. 

\item Further, we provide a similar high-dimensional characterisation for the performance of the optimal Bayesian estimator in this problem. These two results follow from an extension of the replica method for generalised linear estimation to covariate distributions in the family of Eq.~\eqref{eq:def:superstats}, and hold for a broader mixture data model which we discuss in Section~\ref{sec:gencluster}. 

\item Leveraging the characterisation above, we investigate the impact of heavy-tailed covariate and response contamination on the performance of M-estimators, including, for the first time for this model, analytic control of infinite-variance covariates. In particular, our result highlights the necessity of regularising the Huber loss in the high-dimensional regime to achieve optimal performance, as we show that it can be suboptimal in the presence of heavy-tailed response contamination even when the location parameter is optimally tuned, at odds with the optimality results of \citet{Huber64} and \citet{Hampel1974} in the classical regime where $n\gg d$.

\item We show that, despite the strong impact of heavy-tailed contamination in the high-dimensional regime where $\sfrac{n}{d} = \Theta(1)$, the error decay rates $\|\hat{\bbeta}-\bbeta_{\star}\|_{2} = \Theta(n^{-1/2})$ of optimally regularised Huber and least-squares regression are optimal, provided the second moments of both the covariates and the noise are bounded. The results also include an exact expression for the coefficients for the leading order of the decay error rates. In contrast, we show that when the covariates' second moment does not exist (but under the assumption of a finite second moment for the noise), the rates explicitly depend on the tail behavior of the covariates' distribution.
\end{itemize}

\paragraph{Related works ---} Robust regression is a classical topic in statistics, with several books dedicated to the subject \citep{huber2004robust, hampel2011robust, rousseeuw2005robust, maronna2019robust}.  In contrast to the classical regime, literature on robust regression in the high-dimensional regime remains relatively scarce. Early works in this direction provided a characterisation of M-estimators in the proportional limit of $n,d\to\infty$ with fixed $\alpha = \sfrac{n}{d}$ for Gaussian \citep{Karoui2013a, Karoui2013b, Donoho2016} and sub-Gaussian \citep{Karoui2013b} designs with bounded noise variance. In particular, it was shown that in this regime the maximum likelihood estimator is not necessarily the optimal choice of objective $\rho$ \citep{Karoui2013b}. De-biasing and confidence intervals for high-dimensional M-estimation on Gaussian covariates were discussed by \cite{Bellec2022b, Bellec2023b}. The family of distributions defined in Eq.~\eqref{eq:def:superstats} has previously appeared under different names and literatures, such as \emph{superstatistics} in the context of statistical physics \citep{Beck2003}, \emph{elliptical distributions} in the context of statistics \citep{Couillet2015, Karoui2018b, adomaityte2023} and \emph{Gaussian scale mixture} in the context of signal processing \citep{Wainwright1999}. It has been studied in the context of high-dimensional robust covariance estimation by \citet{Couillet2015, Couillet2016}. \citet{Karoui2018b} considered (unregularised) robust regression for elliptically distributed covariates under boundedness conditions on the moments. Our work differs in many directions: (a) we relax the assumption of bounded moments and extend the analysis to the more general case of a \textit{mixture} of elliptical distributions; (b) we derive the asymptotic performance of the \emph{Bayes optimal estimator} for this family;  (c) we consider a generic convex penalty; (d) we investigate the impact of Huber $\epsilon$-contamination.  
More recently, \citet{Vilucchio2023} studied the model Eq.~\eqref{eq:def:data} under a Gaussian design and a double Gaussian noise model for outliers in the proportional high-dimensional regime, showing that estimators of the type in Eq.~\eqref{eq:def:erm} can fail to be consistent. Beyond (sub-)Gaussian designs, high-dimensional upper bounds on the regressor mean-squared error were obtained under different settings \citep{Hsu2016, Lugosi2019, Roy2021}, including heavy-tailed noise \citep{Sun2020} and heavy-tailed covariate contamination \citep{sasai2022robust, Pensia2021}.

%% file: sections/asymptotics.tex
In this section, we discuss our two main theoretical results: the high-dimensional asymptotic characterisation of the M-estimator defined in Eq.~\eqref{eq:def:erm} and the corresponding Bayes-optimal error. In fact, the results in this section hold under the following slightly more general assumptions.

\begin{assumption}[Data] 
\label{ass:data}
The covariates $\bx_{i}\in\R^{d}$, $i\in[n]$ are independently drawn from the family of ``superstatistical'' distributions $\bx_i\sim p(\bx)\coloneqq\mathbb E[\mathcal N(\bx;\mathbf 0,\sfrac{\sigma^2}{d}\bI_d)]$, where the expectation is over $\sigma$ with generic distribution density $\varrho$ supported on the positive real line $\R^*_+\coloneqq (0,+\infty)$. For each $i\in[n]$, the corresponding response $y_{i}\in\mathcal{Y}$ is drawn from a conditional law $P_{0}$ on $\mathcal{Y}$ as in Eq.~\eqref{eq:labels}, with target weights $\bbeta_{\star}\in\R^{d}$ having finite normalised norm $\beta_\star^2\coloneqq \lim_{d\to\infty}\sfrac{1}{d}\|\bbeta_\star\|_2^2$.
\end{assumption}
\noindent Note that the linear model in Eq.~\eqref{eq:def:data} corresponds to the choice $P_0(y|\tau)=p_\eta(y-\tau)$, $p_\eta$ being the density of the noise $\eta$.
\begin{assumption}[Predictor] 
\label{ass:predictor}
We consider the hypothesis class of generalised linear predictors $\mathcal{H} = \{f_{\bbeta}(\bx)\coloneqq f(\bbeta^\intercal\bx), \bbeta\in\R^{d}\}$, where $f\colon\R\to\mathcal Y$ is a generic activation function, and the weights $\bbeta\in\R^{d}$ are obtained by minimising the following empirical risk: 
\begin{equation}
\hat{\bbeta}_{\lambda}\coloneqq \underset{\bbeta\in \R^{d}}{\argmin}~ \sum\limits_{i=1}^n\rho\left(y_i-\bbeta^\intercal\bx_i\right) + \frac{\lambda}{2} \|\bbeta\|_2^2,\qquad \lambda\in\R_+,
\end{equation}
for a convex objective function $\rho\colon\R\to \R_+$.
\end{assumption}

\begin{assumption}[Proportional regime] 
\label{ass:regime}
We consider the proportional high-dimensional regime where both $n,d\to\infty$ at a fixed ratio $\alpha \coloneqq \sfrac{n}{d}$, known as the sample complexity.
\end{assumption}
In particular, note that Assumption \ref{ass:data} covers the additive noise model in Eq.~\eqref{eq:def:data}, and that the standard robust regression setting is given by taking $f(\bbeta^{\intercal}\bx) = \bbeta^{\intercal}\bx$. In Appendix \ref{app:replica} and Appendix \ref{app:ridge_regul} we derive the following result.

\begin{result}[High-dimensional asymptotics]\label{res:K1} 
Let $\varphi\colon\mathcal Y\times \R\to\R$ denote a test function, and define the following generalisation and training statistics:
\begin{equation}\label{eq:errdef}\textstyle
\mathcal{E}_{g}(\hat{\bbeta}) = \mathbb{E}_{(y,\bx)}\big[\varphi(y,\hat\bbeta_{\lambda}^\intercal\bx)\big], \qquad \mathcal{E}_{t}(\hat{\bbeta})\coloneqq \frac{1}{n}\sum\limits_{i=1}^n\varphi(y_i,\hat\bbeta_{\lambda}^\intercal \bx_i).
\end{equation}
Then, under Assumptions \ref{ass:data}--\ref{ass:regime} we have:
\begin{equation}\label{eq:limiterr}
    \mathcal{E}_{g}(\hat{\bbeta})  \xrightarrow[n,d\to\infty]{\rm P} \varepsilon_g(\alpha,\lambda,\beta_{\star}^{2}), \qquad \mathcal{E}_{t} (\hat{\bbeta})\xrightarrow[n,d\to\infty]{\rm P} \varepsilon_t(\alpha,\lambda,\beta_{\star}^{2}),
\end{equation}
if the limit exists. The expressions above are explicitly given by:
\begin{equation}
\begin{split}
\varepsilon_t(\alpha,\lambda,\beta_{\star}^{2})&\textstyle=\int_{\mathcal Y}\dd y\,\mathbb E_{\sigma,\zeta}\left[Z_0\left(y,\frac{m\sigma}{\sqrt{q}}\zeta,\sigma^2\beta_\star^2-\frac{\sigma^2m^2}{q}\right)\varphi\left(y,\sigma\sqrt{q}\zeta+v\sigma^2f\right)\right],\\\
\varepsilon_g(\alpha,\lambda,\beta_{\star}^{2})&\textstyle=\int_{\mathcal Y}\dd y\int\dd\eta\int\dd\tau\,P_0\left(y|\tau\right)\mathbb E_{\sigma}\left[\mathcal N\left(\begin{psmallmatrix}{\tau}\\{\eta}\end{psmallmatrix};\mathbf 0,\sigma^2\begin{psmallmatrix}\beta_\star^2&m\\m&q\end{psmallmatrix}\right)\right]\varphi(y,\eta),
\end{split}\end{equation}
where $\bxi\sim\mathcal N(\mathbf 0,\bI_d)$, $\zeta\sim\mathcal N(0,1)$ are independent Gaussian variables. Similarly, the mean-squared error on the estimator is given by:
\begin{equation}\comprimi
\label{eq:def:esterr}
\varepsilon_{\rm est}(\alpha,\lambda,\beta_{\star}^{2})\coloneqq\lim_{n,d\to+\infty}\frac{1}{d}\mathbb E_{\mathcal D}\big[\|\hat\bbeta_\lambda-\bbeta_\star\|_2^2\big]=\beta_\star^2-2m+q\xrightarrow{\alpha\to+\infty}\beta_\star^{-2}\lim_{\alpha\to+\infty}(m-\beta_\star^2)^2.
\end{equation}
Moreover, the angle between the estimator $\hat\bbeta_\lambda$ and $\bbeta_\star$ is given by
\begin{equation}\textstyle
\lim_{n,d\to+\infty}\mathbb E_{\mathcal D}[\mathrm{angle}(\hat\bbeta_\lambda,\bbeta_\star)]=\frac{1}{\pi}\arccos\left(\frac{m}{\beta_\star\sqrt{q}}\right)\xrightarrow{\alpha\to+\infty}0.
\end{equation}
In all formulas above we have introduced the order parameters $v$, $q$, and $m$. These quantities are found by solving the following set of self-consistent fixed-point equations:
\begin{subequations}\label{eq:sp}
\begin{equation}
\begin{split}
	q &=\textstyle\frac{\hat m^2\beta_\star^2+\hat q}{(\lambda+\hat v)^2}\\
	m &=\textstyle \frac{\beta_\star^2\hat m}{\lambda+\hat v}\\
	v &=\textstyle\frac{1}{\lambda+\hat v}
\end{split}\qquad
\begin{split}
\hat q&=\textstyle\alpha \int_{\mathcal Y}\dd y\,\mathbb E_{\sigma,\zeta}\left[\sigma^2 Z_0\left(y,\frac{m\sigma}{\sqrt{q}}\zeta,\sigma^2\beta_\star^2-\frac{\sigma^2m^2}{q}\right)f^2\right],\\
\hat v&=\textstyle-\alpha\int_{\mathcal Y}\dd y\,\mathbb E_{\sigma,\zeta}\left[\sigma^2 Z_0\left(y,\frac{m\sigma}{\sqrt{q}}\zeta,\sigma^2\beta_\star^2-\frac{\sigma^2m^2}{q}\right)\partial_\omega f\right],\\
\hat m&=\textstyle\alpha \int_{\mathcal Y}\dd y\,\mathbb E_{\sigma,\zeta}\left[\sigma^2\partial_\mu Z_0\left(y,\frac{m\sigma}{\sqrt{q}}\zeta,\sigma^2\beta_\star^2-\frac{\sigma^2m^2}{q}\right) f\right],
\end{split}   
\end{equation}
where the auxiliary function $Z_{0}$ and proximal operator $f$ are defined by:
\begin{equation}
\label{eq:def:Z0}
\begin{split}
Z_0(y,\mu,v)&\coloneqq\mathbb E_z[P_0(y|\mu+\sqrt{v}z)],\\
f&\textstyle\coloneqq\arg\min_{u}\left[\frac{\sigma^2vu^2}{2}+\rho\left(y-\sigma\sqrt{q}\zeta-\sigma^2 v u\right)\right]\in\R.
\end{split}
\end{equation}
\end{subequations}
\end{result}
Result \ref{res:K1} reduces the high-dimensional optimisation problem defined by Eq.~\eqref{ass:predictor} to a set of low dimensional self-consistent equations for the parameters $\hat{q},q,\hat{m},m,\hat{v}, v\in\R$. Despite being cumbersome, this low-dimensional problem can be efficiently solved numerically. Note that our result for $\varepsilon_{\rm est}$ implies that the estimator $\hat\bbeta_\lambda$ is consistent if and only if $\lim_{\alpha\to+\infty}m=\beta_\star^2$. In Appendix \ref{app:square} we explicitly derive the convergence rates of $\varepsilon_{\rm est}$ for the cases of square loss and {Huber loss}, obtaining the following.

\begin{result}[Convergence rates] \label{res:scaling}
Consider the linear model in Eq.~\eqref{eq:def:data} with covariates as in Eq.~\eqref{eq:def:superstats}. Assume that $\mathbb{E}[\eta_i]=0$, $\hsigma_0^2 \coloneqq \mathbb{E}[\eta_i^2]\in\R^*_+$. Let us also assume that in Eq.~\eqref{eq:def:superstats} we have $\varrho(\sigma)\sim\frac{1}{\sigma^{2a+1}}$ for $\sigma\gg 1$ with $a>0$. 
Then, the ridge estimator $\hat\bbeta_{\lambda}$ minimising Eq.~\eqref{eq:def:erm} with $\rho(t)=\frac{1}{2}t^{2}$ and the estimator obtained using the Huber loss in Eq.~\eqref{eq:def:huberloss} are consistent and satisfy the following scaling:
\begin{equation}\label{eq:result25}\comprimi
\varepsilon_{\rm est}\underset{\alpha \gg 1}{\sim}\begin{cases}
\frac{1}{\alpha}+o\left(\frac{1}{\alpha}\right)&\text{if }a>1,\\
\frac{1}{\alpha\ln\alpha}+o\left(\frac{1}{\alpha\ln\alpha}\right)&\text{if }a=1,\\
\frac{1}{\alpha^{\sfrac{1}{a}}}+o\left(\frac{1}{\alpha^{\sfrac{1}{a}}}\right)&\text{if }a\in(0,1).
\end{cases}
\end{equation}
In particular, in the case of the square loss, the coefficients of the scaling can be explicitly given in terms of the Stieltjes transform $\mathsf S_{\sigma^2}(x)\coloneqq\mathbb E\big[\frac{1}{x+\sigma^2}\big]$ of the measure of $\sigma^2$ as
\begin{equation}\label{eq:result25b}\comprimi
\varepsilon_{\rm est}\underset{\alpha \gg 1}{=}\begin{cases}
\frac{\hsigma_0^2}{\sigma^2_0\alpha}+o\left(\frac{1}{\alpha}\right)&\text{if }a>1,\quad\text{with }\sigma_0^2\coloneqq\mathbb E[\sigma^2]=\lim\limits_{x\to +\infty}x\left(1-x\mathsf S_{\sigma^2}(x)\right),\\
\frac{\hsigma_0^2}{\tsigma^2_0\alpha\ln\alpha}+o\left(\frac{1}{\alpha\ln\alpha}\right)&\text{if }a=1,\quad\text{with }\tsigma_0^2=\lim\limits_{x\to +\infty}\frac{x}{\ln x}\left(1-x\mathsf S_{\sigma^2}(x)\right),\\
\frac{\hsigma_0^2}{(\tsigma_0^2\alpha)^{\sfrac{1}{a}}}+o\left(\frac{1}{\alpha^{\sfrac{1}{a}}}\right)&\text{if }a\in(0,1),\quad\text{with }\tsigma_0^2=\lim\limits_{x\to +\infty}x^{a}\left(1-x\mathsf S_{\sigma^2}(x)\right).
\end{cases}
\end{equation}
\end{result}

By definition, the mean-squared error $\varepsilon_{\rm est}$ of the M-estimator in Eq.~\eqref{eq:def:esterr} is lower-bounded by the minimum mean-squared error achieved by the Bayes-optimal estimator $\hat\bbeta$, which is given by the expectation over the Bayesian posterior $\hat\bbeta=\mathbb E[\bbeta|\mathcal D]$. Our second main contribution in this manuscript, derived in Appendix \ref{app:BO}, is to provide a sharp asymptotic formula for the minimum mean-squared error under the assumption that $\bbeta_\star\sim \mathcal{N}(\bx;\mathbf 0,\beta_{\star}^{2}\bI_d)$.
\begin{result}[Bayes optimal error]
\label{res:BO} 
Consider the minimum mean-squared error for the data model defined by \ref{ass:data} with $\bbeta_\star\sim \mathcal{N}(\bx;\mathbf 0,\beta_{\star}^{2}\bI_d)$:
\begin{equation}\textstyle
    \mathcal{E}_{\rm BO} \coloneqq \underset{\bbeta\in\R^{d}}{\min}~\frac{1}{d}\mathbb{E}[\|\bbeta_{\star}-\bbeta\|_2^{2}].
\end{equation}
Then, in the high-dimensional limit defined by Assumption \ref{ass:regime}:
\begin{equation}\label{eq:esterr_BO}\textstyle
\mathcal{E}_{\rm BO} \xrightarrow[n,d\to\infty]{\rm P}\varepsilon_{\rm BO} = \beta_\star^2-\sq,
\end{equation}
where $\sq$ is given by the solution of the following self-consistent equations:
\begin{equation}
\shq=\alpha  \int_{\mathcal Y}\dd y\,\mathbb E_{\sigma,\zeta}\Big[\sigma^2 Z_0(y,\mu,V)\left(\partial_\mu\ln Z_0(y,\mu,V)\right)^2\Big|_{\substack{\mu=\sigma\sqrt{\sq}\zeta\\ V=\sigma^2(\beta_\star^2-\sq)}}\Big],\quad 
	\sq =\frac{\beta_\star^4\shq}{1+\beta_\star^2\shq}
\end{equation}
with $\zeta\sim\mathcal N(0,1)$ and $Z_{0}$ given in Eq.~\eqref{eq:def:Z0}.
\end{result}
Similarly to Result \ref{res:K1}, Result \ref{res:BO} reduces the problem of estimating the Bayes optimal error from the computationally intractable evaluation of the high-dimensional posterior marginals to a simple two-dimensional set of self-consistent equations on the variables $\shq, \sq\in\R$. 

A detailed derivation of both Results \ref{res:K1} \& \ref{res:BO} is provided in Appendix \ref{app:replica} using the replica method. Despite being non-rigorous, Results \ref{res:K1} \& \ref{res:BO} provide a natural extension of rigorous results in the established literature of high-dimensional asymptotics for generalised linear estimation on Gaussian covariates \citep{Karoui2013b, Donoho2016, Thrampoulidis2018, Barbier2019, Loureiro_2022} to the elliptical family defined in Assumption~\ref{ass:data} . Indeed, by taking $\varrho(\sigma) = \delta_{\sigma,\bar\sigma}$ for some $\bar\sigma\in\R^*_+$, Result \ref{res:K1} reduces to the rigorous formulas proven by \cite{Thrampoulidis2018}, and Result \ref{res:BO} reduces to the rigorous formulas proven by \cite{Barbier2019}. Nevertheless, as we will see in the next section, they produce excellent agreement with moderately finite size simulations, across different choices of penalty $\rho$ and covariate distribution, including cases in which the variance is infinite, see Fig. \ref{sec3:fig:label_noise_contamination} for instance.

%% file: sections/discussion.tex
In this section, we investigate the consequences of Results \ref{res:K1} \& \ref{res:BO} in the context of robust regression with heavy-tail contamination as in Eq.~\eqref{eq:def:hubercont} within the model in Eq.~\eqref{eq:def:data}. Our theoretical results are given for \textit{any} noise distribution $p_\eta$, but, to exemplify our findings, we will focus on the case of heavy-tailed contamination of standard Gaussian noise, in a form similar to Eq.~\eqref{eq:def:hubercont}: 
\begin{align}
\label{eq:def:hubercontnoise}
p_\eta(\eta)&=\mathbb E_\sigma[\mathcal N(\eta;0,\sigma^2)], \\\text{with}\quad
\sigma&\sim\hat\varrho=(1-\epsilon_{\rm n})\delta_{\sigma,1}+\epsilon_{\rm n}\hat\varrho_0,\quad\epsilon_{\rm n}\in[0,1] \nonumber.   
\end{align}
As in the case of covariates, $\epsilon_{\rm n}$ measures the level of contamination of the purely Gaussian noise.  For concreteness, we take $\bbeta_{\star}\sim\mathcal{N}(0,\beta_{\star}^{2}\bI_{d})$ and use for $\varrho_0$ and $\hat\varrho_0$ one of the distributions in Table \ref{tab:examples}, which provide convenient parametric families of distributions where the existence of moments is easily tunable.

\begin{wrapfigure}{r}{0.4\textwidth}
\vspace{-5mm}\includegraphics[width=0.4\textwidth]{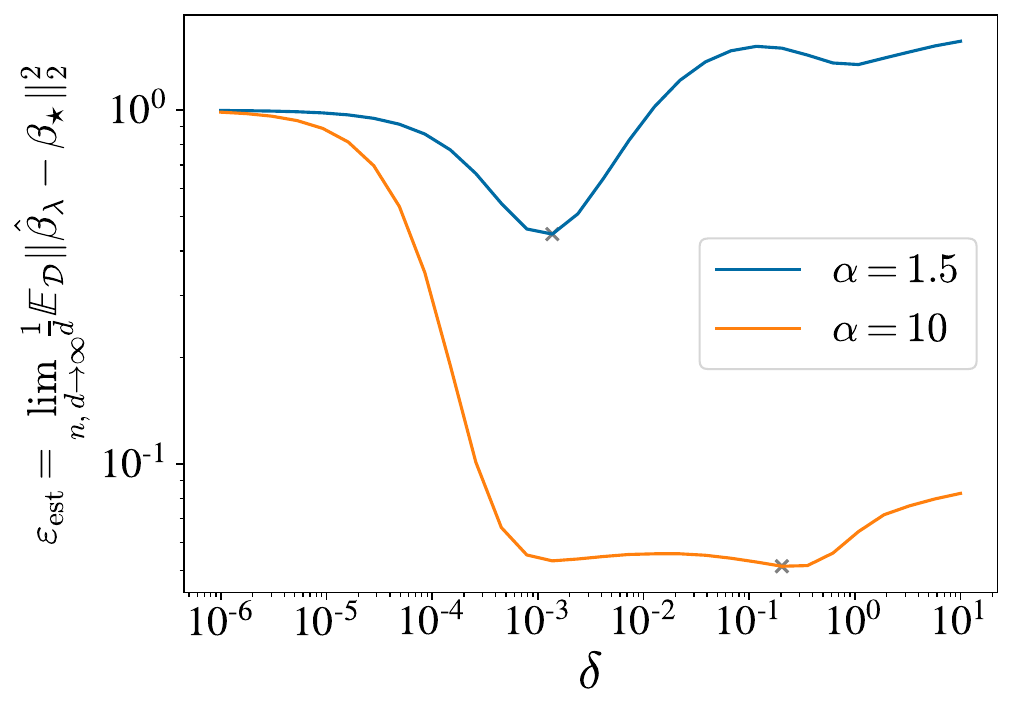}
\caption{Value of $\varepsilon_{\rm est}$ obtained using a regularised Huber at given $\lambda=10^{-3}$ as a function of $\delta$ for different values of $\alpha$. Here the contamination level is $\epsilon_{\rm n}=0.5$.}
\label{sec3:fig:mindelta}
\end{wrapfigure}

\subsection{Heavy-tailed label contamination}
\label{sec:noise}
In this section, we discuss the impact of \textit{label contamination} in high-dimensional robust regression. Since the Gaussian assumption for the noise label is widespread in statistics, we will focus our discussion on an $\epsilon$-contamination of the standard normal distribution. We do this by assuming a noise distribution $p_\eta$ as in Eq.~\eqref{eq:def:hubercontnoise} and adopting a suitable contaminating distribution $\hat\varrho_0$ that generates heavy-tailed label noise when $\epsilon_{\rm c}>0$. We recall the reader that when the loss function is the negative log-likelihood of the noise $-\log p_\eta$, strong guarantees hold for the M-estimator in the classical limit $n\gg d$. In particular, the estimator is unbiased and asymptotically normal, with estimation error (defining the variance)  $\|\bbeta_{\star}-\hat\bbeta_\lambda\|_2^2=\Theta(n^{-1})$ \citep{van2000asymptotic}. Since the square loss corresponds to the negative log-likelihood of standard Gaussian noise, in the absence of contamination $\epsilon_{\rm n}=0$ the square loss attains an optimal rate as $n\to\infty$. Our exact high-dimensional characterisation in terms of Results \ref{res:K1} \& \ref{res:BO} allows to check whether the square loss remains optimal in high dimensions, i.e., when $n=\Theta(d)$, assuming contaminated labels as in Eq.~\eqref{eq:def:hubercontnoise}. On the other hand, we will assume Gaussian covariates in this section, i.e., $\epsilon_{\rm c}=0$ in Eq.~\eqref{eq:def:hubercont}. 

\begin{figure}[t]
    \centering
    \includegraphics[width=\textwidth]{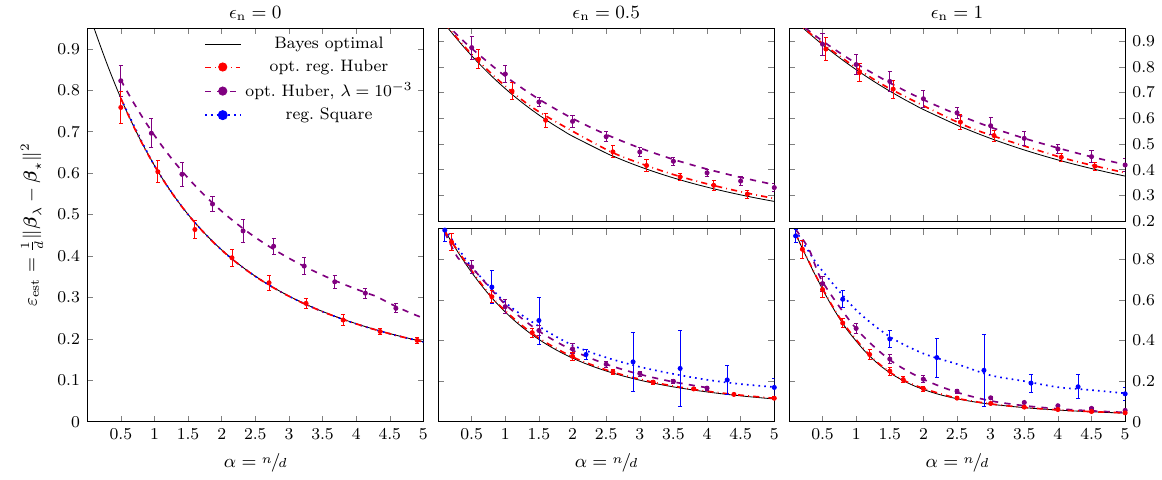}
    \caption{
    Label contamination of standard Gaussian covariates, with distribution $p(\bx)=\mathcal N(\bx;0,\sfrac{1}{d} \bI_d)$, as a function of the sample complexity $\alpha=\sfrac{n}{d}$ for different contamination levels $\epsilon_{\rm n}\in[0,1]$ as in Eq.~\eqref{eq:def:hubercontnoise}. The label-contaminating distribution is inverse-gamma $\hat\varrho_0(\sigma)\propto\sigma^{-2a-1}\exp(-b\sigma^{-2})$, for details see Table~\ref{tab:examples}. Theoretical predictions (lines) are compared with the results of numerical experiments (dots) obtained averaging over $20$ instances with $d=10^3$. (\textbf{Left}) Purely Gaussian noise $(\epsilon_{\rm n}=0)$. (\textbf{Top}). Case $\epsilon_{\rm n}>0$, $b=1$ and $a=\sfrac{4}{5}<1$, implying $\hat\varrho_0(\sigma)$ has infinite mean and thus $\mathbb E[\eta^2]=+\infty$ but $\mathbb E[\eta]<+\infty$. The performance degrades as the contamination level $\epsilon_{\rm n}$ is increased. Optimally regularised Huber (red) achieves the optimal Bayesian performance (solid), while with fixed small regularisation the Huber estimator performs suboptimally  (purple). Square loss results are not represented as, in this case, the average estimation error is not finite. (\textbf{Bottom}). Case $\epsilon_{\rm n}>0$, $a=1+b=1+\sfrac{1}{10}$, corresponding to $\mathbb E[\eta^2]=1$. The performance uniformly improves as the contamination $\epsilon_{\rm n}$ grows. Optimally regularised Huber (red) achieves optimal Bayesian performance (solid), while both Huber with untuned regularisation (purple) and optimally regularised ridge (blue) are suboptimal.
    }\label{sec3:fig:label_noise_contamination}
\end{figure}
\paragraph{The importance of regularisation ---} 
Fig.~\ref{sec3:fig:label_noise_contamination} shows the performance of the ridge and Huber estimators in terms of the estimation error $\varepsilon_{\rm est}$ in Eq.~\eqref{eq:def:esterr}, and the corresponding estimation error $\varepsilon_{\rm BO}$ of the Bayes optimal estimator in Eq.~\eqref{eq:esterr_BO}, for various values of the sample complexity $\alpha=\sfrac{n}{d}$ which can be intuitively understood as a signal-to-noise ratio. The theoretical lines and experiments (dots) show good agreement. The left panel in Fig.~\ref{sec3:fig:label_noise_contamination} shows the performance when the noise is purely Gaussian, i.e., $\epsilon_{\rm n}=0$ in  Eq.~\eqref{eq:def:hubercontnoise}. As expected, in the absence of contamination, optimally regularised ridge regression achieves the minimum mean-squared error $\varepsilon_{\rm est}$ for this model. Indeed, in this case, the square loss is not only the maximum likelihood estimator but also coincides with the Bayes-optimal estimator, since from a Bayesian perspective the $\ell_2$ penalty corresponds to the optimal Bayesian prior when optimally tuned. The center and right panels in Fig.~\ref{sec3:fig:label_noise_contamination} illustrate two distinct scenarios, i.e., contamination of the labels with an infinite-variance noise (top), so that $\mathbb E[\eta^2]=+\infty$, and contamination with a finite-variance noise (bottom), so that $\mathbb E[\eta^2]=1$. While in the former case the performance degrades as a function of the contamination, counter-intuitively, it improves when the contamination variance is bounded. As contamination is introduced $(\epsilon_{\rm n}>0)$, optimally regularised ridge regression is observed to be suboptimal in the finite-variance scenario. More surprising, perhaps, is that, at small values of the regularisation $\lambda$, the Huber loss in Eq.~\eqref{eq:def:huberloss} with optimally chosen location parameter $\delta^{\star}$ is also suboptimal in the high-dimensional regime $n=\Theta(d)$. This is to be contrasted with the well-known optimality results in the classical regime $n\gg d$, which are recovered at $\alpha\to\infty$. Interestingly, at small but finite $\lambda$, the improvement of performance for the Huber loss with $\alpha$ coincides with a sharp transition of the optimal location parameter $\delta^*$ from $O(\lambda)$ to $O(1)$ at a given value of $\alpha$: the value of $\varepsilon_{\rm est}$ develops indeed two minima as a function of $\delta$, one in $\delta=O(\lambda)$ and another in $\delta=O(1)$, whose relative depth changes with $\alpha$, see Fig.~\ref{sec3:fig:mindelta}. For $\lambda\to 0$, the minimum in $\delta=O(\lambda)$ disappears. This phenomenology is observed for various choices of $\hat\varrho_0$ (see Appendix \ref{app:furthernumerics}) but it also persists in other settings, as we mention later in the text. Additionally, taking $\delta\rightarrow 0$ in Eq.~\eqref{eq:def:huberloss} corresponds to the $\ell_1$ loss $\rho(t)=|t|$ in Eq.~\eqref{eq:def:erm}, leading to the least absolute deviation (LAD) estimator, which performs similar or worse than other estimators considered, and is discussed in Appendix~\ref{app:sec:LAD}. In general, adding a $\ell_{2}$ regularisation and cross-validating significantly improves the performance of the Huber loss, bringing it to or very close to Bayes optimality at all contamination levels. Our results highlight the necessity of properly regularising robust estimators in order to achieve optimality in the high-dimensional regime. In particular, optimally regularising Huber is crucial to achieve Bayes-optimal performances in the case of infinite-variance label noise. 

\begin{wrapfigure}[20]{L}{0.4\textwidth}
    \centering \vspace{-5mm}
    \includegraphics[width=0.4\columnwidth]{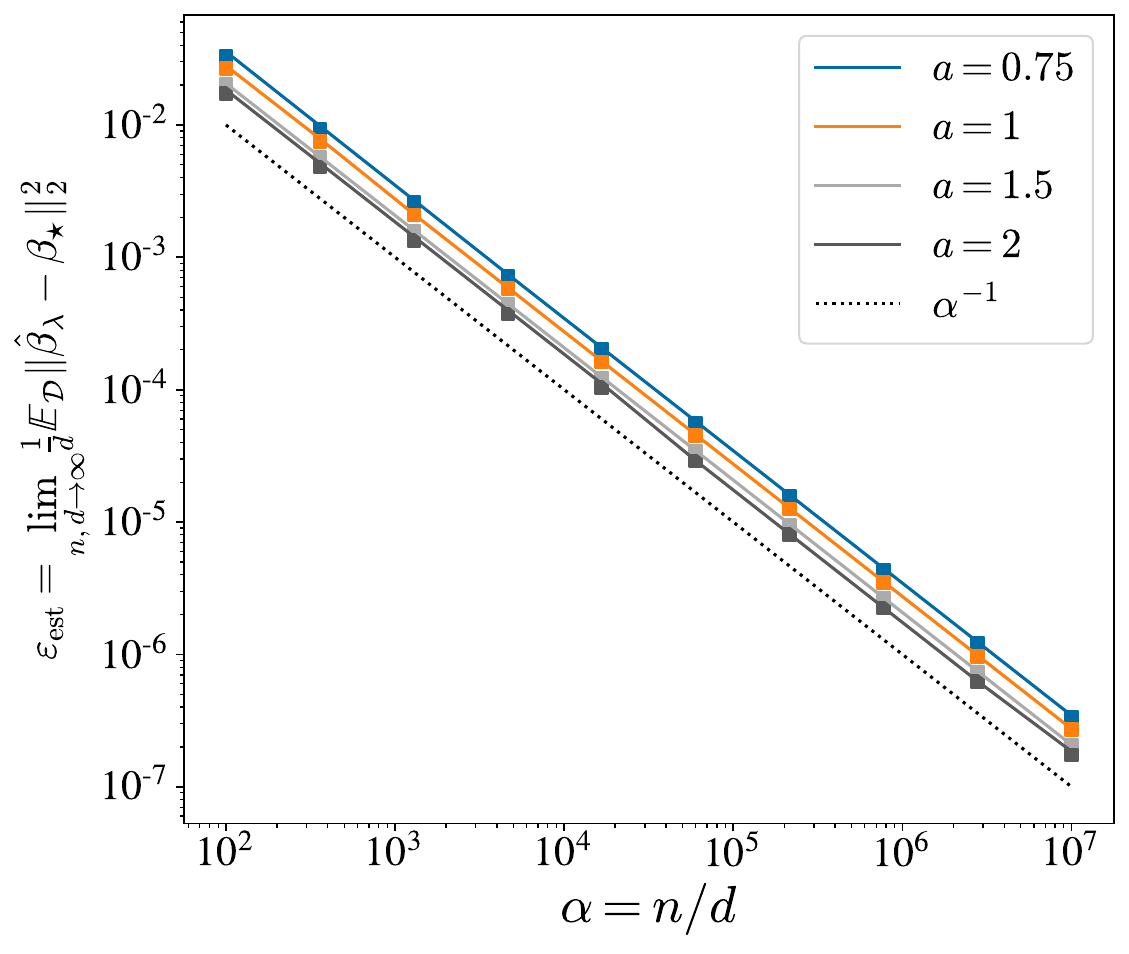}
    \caption{Estimation error $\varepsilon_{\rm est}$  as in Eq.~\eqref{eq:def:esterr} at large-$\alpha$ using regularised optimal Huber loss (solid lines). The covariates are Gaussian, whereas the label noise is obtained as in Eq.~\eqref{eq:def:hubercontnoise} with $\varrho(\sigma)\sim\sigma^{-2a-1}$,  $a>0$. The results are compared with the Bayes-optimal performance (squares). The dotted line shows a scaling of $\alpha^{-1}$.}\label{sec3:fig:result25_check}\vspace{-2cm}
\end{wrapfigure}

\paragraph{Convergence rates ---} In Fig.~\ref{sec3:fig:result25_check} we consider the case of full contamination $\epsilon_{\rm n}=1$ of the labels with both finite- and infinite-variance noise in the large-$\alpha$ regime: this corresponds to the classical limit $n \gg d$. In particular, we focused on the regularised Huber loss, whose location parameter and regularisation strength were optimised to reach the minimum estimation error. A heavier tail in the noise hurts performance in terms of the estimation error $\varepsilon_{\rm est}$, but does not affect its $\alpha^{-1}$ scaling, which remains universal even when $\mathbb E[\eta^2]=+\infty$ and is observed also in the Bayes optimal lower bound. We will see below that this will not be the case if a contamination of the covariates is considered. 

\subsection{Covariate contamination in the presence of heavy-tailed label noise} \label{sec:covariatecont}
\begin{figure}[t]
    \centering
    \includegraphics[height=0.355\textwidth]{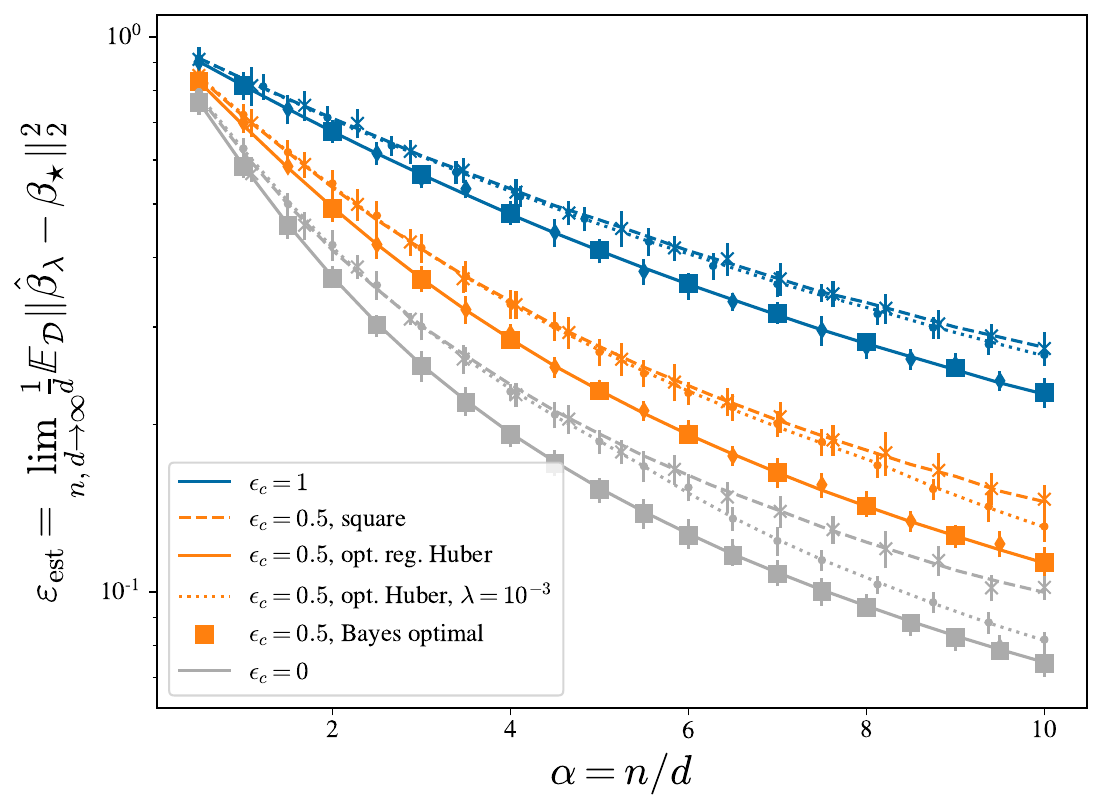}
    \quad
    \includegraphics[height=0.355\textwidth]{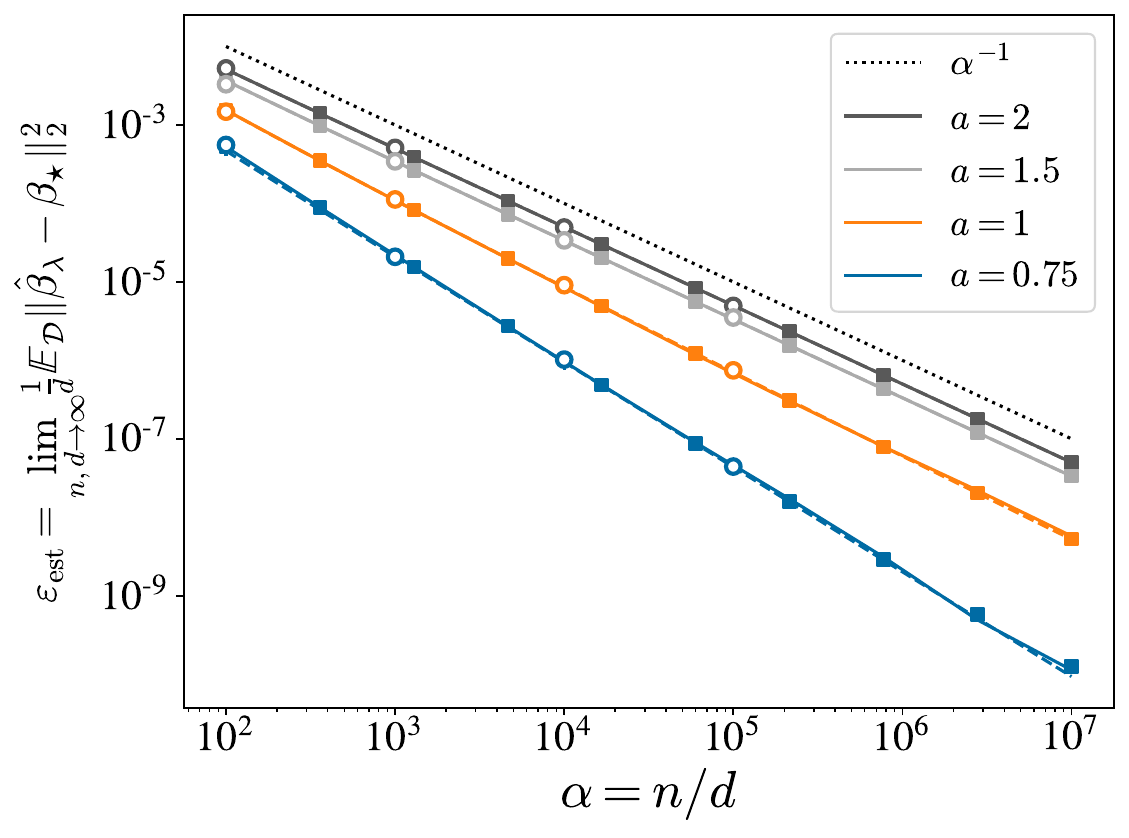}
    \caption{(\textbf{Left}) Estimation error as a function of sample complexity $\alpha=\sfrac{n}{d}$. Dots correspond to the average test error of $50$ numerical experiments in dimension $d=10^3$. Covariates are contaminated as in Eq.~\eqref{eq:def:hubercont} using for the contaminating distribution $\varrho_0$ an inverse-gamma distribution with $b=a-1=\sfrac{1}{10}$, see Table~\ref{tab:examples} for details. The labels are distributed according to an inverse-gamma distribution ($\epsilon_{\rm n}=1$) with parameters $a=1+b=2$, corresponding to $\mathbb E[\eta^2]=1$. (\textbf{Right}) Estimation error $\varepsilon_{\rm est}$ as in Eq.~\eqref{eq:def:esterr} at large-$\alpha$ obtained from our theory for the square loss (dashed line), regularised optimal Huber loss (solid line) and Bayes-optimal performance given by Result~\ref{res:BO} (squares). The covariates' variance is Pareto-distributed to have $\varrho(\sigma)\sim\sigma^{-2a-1}$ with $a>0$, and the label noise is Gaussian. White dots correspond to numerical experiments in dimension $d=50$, averaged over $50$ instances. The black dotted line shows a scaling of $\alpha^{-1}$. 
    }\label{sec3:fig:cov_cont_rates}
\end{figure}
We now move to discuss the impact of \emph{covariates' contamination}, focusing on the case in which covariates are generated with a density as in Eq.~\eqref{eq:def:hubercont}. On top of this, we assume a label noise \textit{with finite variance}, $\mathbb E[\eta^2]=1$. Note that, as discussed above, in this high-dimensional regime, regularised ridge regression is not optimal in the presence of heavy-tailed label noise.   

Fig.~\ref{sec3:fig:cov_cont_rates} (left) focuses on a covariate-contaminated model with Gaussian noise. It shows that as the contamination level of the covariates $\epsilon_{\rm c}$ grows, the performance is negatively affected for all metrics. The phenomenology is similar to the one observed for label contamination: optimally-tuned Huber loss achieves Bayes-optimal performance for optimal regularisation, while optimally regularised square loss performs worst. An optimally tuned Huber loss ridge-regularised with a given strength $\lambda$ performs halfway between the square and the Bayes-optimal bound, approaching the latter as $\alpha$ grows. As observed in the previous section, the improvement in the $\delta$-tuned Huber loss occurs, for small but finite $\lambda$, with a sharp jump in the optimal $\delta^*$ parameter from a value $O(\lambda)$ to a value $O(1)$ as $\alpha$ increases, reflected in a kink in the curve of $\varepsilon_{\rm est}$ as function of $\alpha$. Finally, an optimally $(\delta,\lambda)$-tuned regularised Huber loss achieves Bayes-optimal performance. The phenomenology is analogous to the one described in Section \ref{sec:noise}.

\paragraph{Convergence rates ---} As we show in Result~\ref{res:scaling}, if the statistician adopts the square loss {or the Huber loss} in the study of power-law tailed covariates, the estimation error rate might depend on the tail, under the assumption that the second moment of the noise is finite. In particular, Result \ref{res:K1} implies that, in this case, the estimation error rate depends only on the covariates second moment. As a consequence, the rate $\varepsilon_{\rm est}\sim \alpha^{-1}$ is universal as long as the second moment of the covariates exists. Curiously, if the second moment does not exist, the estimation error decays \emph{faster}, depending explicitly on the tail exponent of the distribution. Remarkably, this suggests that the presence of fat tails in the covariates can actually improve the classical convergence rates of M-estimation when $n\to\infty$. 
We have numerically verified that the convergence rates reported in Result \ref{res:K1} are Bayes-optimal and are also attained by optimally regularising the Huber loss in different scenarios. This is illustrated in Fig.~\ref{sec3:fig:cov_cont_rates} (right), which clearly shows the dependence on the tail exponent of the covariate distribution when the latter has no finite second moment. The theoretical predictions for the square loss (dashed lines) are plotted using the order parameters that are the fixed points of the set of self-consistent equations written in terms of the Stieltjes transform in Eq.~\eqref{app:eq:spC1}, which are in good agreement with numerical experiments in dimension as low as $d=50$.

The careful reader has noted that the case of covariates with no covariance induces a blow-up of the asymptotic risk when using a square loss: this is indeed what our asymptotic formulas predict as well. Nevertheless, the quantities $m$ and $q$ in Result \ref{res:K1}, and therefore the expected estimation error $\varepsilon_{\rm est}$, are finite and well defined. Such theoretical prediction is verified by solving the ERM task via the explicit formula for the estimator in terms of the dataset, avoiding therefore the issue related to the loss fluctuations that appears using, for example, SGD.
\begin{wrapfigure}[27]{t}{0.5\textwidth}\vspace{4mm}
    \includegraphics[width=0.5\columnwidth]{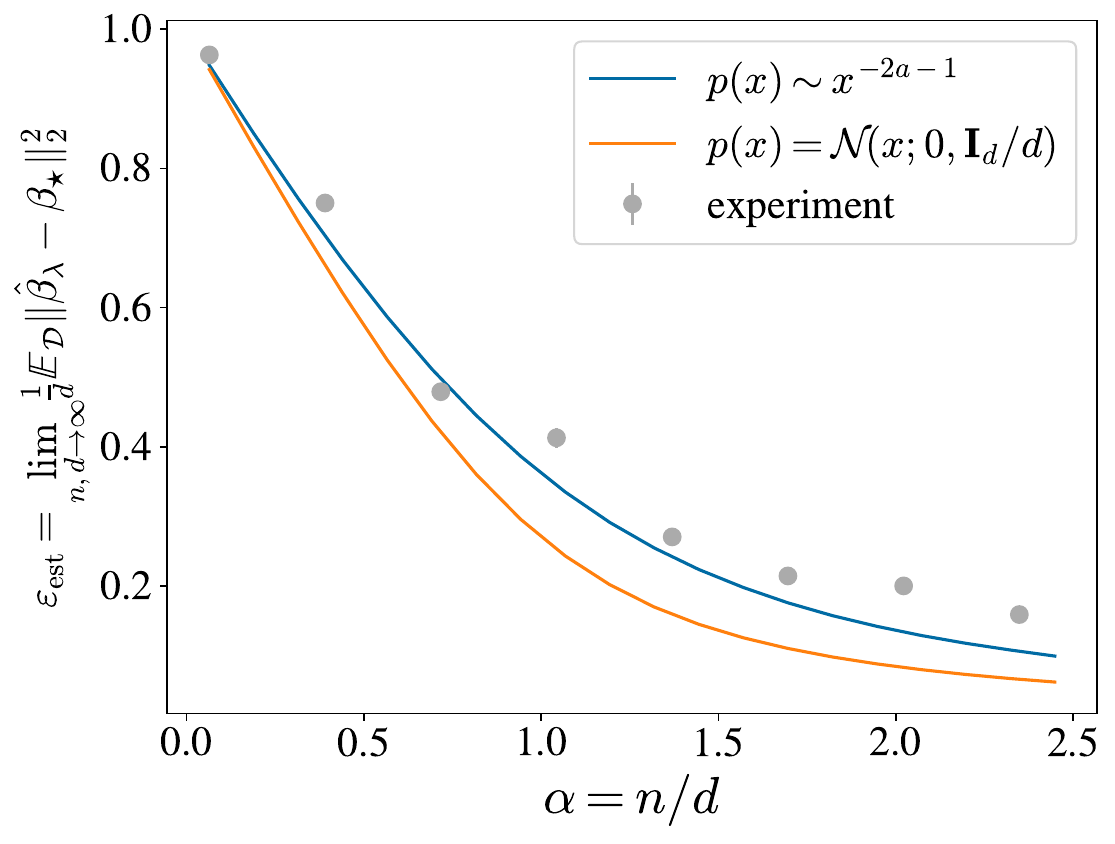}
    \caption{Estimation error $\varepsilon_{\rm est}$ as in Eq.~\eqref{eq:def:esterr} as a function of $\alpha$ using $\ell_2$-regularised square loss with $\lambda=\sfrac{1}{10}$. Dots correspond to the average error of $200$ numerical experiments (varying $n$) on $\bz_i$ with labels generated as Eq.~\eqref{eq:def:data} with Gaussian label noise. Theoretical predictions obtained by taking into account the mean of the covariance spectrum and the covariates' power-law decay (blue) show better agreement with the experiment than a standard Gaussian equivalence assumption (orange).
    }\label{sec3:fig:real_data}
\end{wrapfigure}

\paragraph{An experiment on a real dataset ---} We conclude our numerical tests by comparing our theoretical results with an ERM task on real heavy-tailed data. We considered a time series $\{\bx_i\}_{i=1}^n$ of $d=858$ stocks traded on the NYSE\footnote{The ``Daily Stock Files'' data was accessed from the \citet{CRSP} via WRDS.}: each covariate $\bx_i\in\R^d$ contains the daily returns spanning $n=2519$ days. Stock returns are known to exhibit a power-law decay \cite{Mandelbrot63, Fama65, MantegnaStanley95, Cont_styl_facts, Bouchaud_Potters_2003}. After centering our covariates around their empirical mean, we rescaled them so that the resulting covariance has average eigenvalue $\sfrac{1}{d}$. The newly obtained covariates $\bz_i$ retain a radial power-law decay of the type $p(\bz)\sim\|\bz\|^{-2a-1}$ with $a\simeq1.65$ (see Appendix \ref{app:C3:real_data} for details).
We labeled our dataset by using a linear model as in Eq.~\eqref{eq:def:data}, and adding purely Gaussian label noise $\eta\sim\mathcal{N}(0,\sfrac{1}{10})$. In Fig.~\ref{sec3:fig:real_data} we plot $\varepsilon_{\rm est}$ as in Eq.~\eqref{eq:def:esterr} for $\ell_2$-regularised square loss, comparing the real data results with two theoretical assumptions. First, we assume a `Gaussian equivalent' model by assuming $\bz\sim\mathcal N(\mathbf 0,\sfrac{1}{d}\bI_d)$. Second, we assume heavy-tailed covariates with covariance $\sfrac{1}{d}\bI_d$ but power-law decay produced by selecting inverse-gamma $\varrho$ as in Eq.~\eqref{eq:def:hubercont} with $b=a-1$. 
Despite the rough approximation, the Gaussian theoretical model falls short in the description of the experiment, whilst by simply taking into account the power-law decay of the dataset a much better agreement with the experiments for various $\alpha$ is obtained.

%% file: sections/generalisation.tex
We conclude this contribution by stating the most general form of our exact high-dimensional asymptotic result, valid for an arbitrary mixture of $K$ elliptical distributions, that we call \textit{elliptical mixture model} (EMM). To be more precise, let us consider the following general Assumptions:
\begin{assumption}[EMM Data] 
\label{ass:datassm}
In the dataset $\mathcal D\coloneqq\{(y_i,\bx_i)\}_{i\in[n]}$, the covariates $\bx_{i}\in\R^{d}$, $i\in[n]$ are independently drawn from the mixture $\bx_i\sim p(\bx)\coloneqq\sum_{c=1}^K p_c\mathbb E_{\sigma_c}[\mathcal N(\bx;\bmu_c,\sfrac{\sigma^2_c}{d}\bI_d)]$, where, for each $c\in[K]$, $d\|\bmu_c\|_2^2=\Theta(1)$, $p_c\in[0,1]$ and the expectation is over $\sigma_c\sim\varrho_c$, generic distribution density supported on $\R^*_+$. Moreover, $\sum_{c=1}^Kp_c=1$. For each $\bx_i$, $i\in[n]$, the corresponding response $y_{i}\in\mathcal{Y}$ is drawn from a conditional law $P_{0}(y|\bbeta_\star^\intercal\bx_i)$ on $\mathcal{Y}$, with target weights $\bbeta_{\star}\in\R^{d}$ having finite normalised norm $\beta_\star^2\coloneqq \lim_{d\to\infty}\sfrac{1}{d}\|\bbeta_\star\|_2^2$.
\end{assumption}
\begin{assumption}[Predictor] 
\label{ass:predictorssm}
We consider the hypothesis class of generalised linear predictors $\mathcal{H} = \{f_{\bbeta}(\bx)\coloneqq f(\bbeta^\intercal\bx), \bbeta\in\R^{d}\}$, where $f\colon\R\to\mathcal Y$ is a generic activation function, and the weights $\bbeta\in\R^{d}$ are obtained by minimising the following empirical risk: 
\begin{equation}\textstyle
\hat{\bbeta}_{\lambda}\coloneqq \underset{\bbeta\in \R^{d}}{\argmin}~ \sum\limits_{i=1}^n\ell\left(y_i,\bbeta^\intercal\bx_i\right) + \lambda r(\bbeta),\qquad \lambda\in\R_+,\label{eq:ERMgen}
\end{equation}
for a convex loss function $\ell\colon\mathcal{Y}\times \R\to \R_+$ and a convex regularisation function $r\colon\R\to\R_+$.
\end{assumption}
Under the further Assumption~\ref{ass:regime}, in Appendix \ref{app:replica} it is shown that the following holds.
\begin{result}[High-dimensional asymptotics for the EMM]\label{res:general} Let us consider the independent random variables $\bxi\sim\mathcal N(\textbf{0},\bI_d)$, and $z,\zeta\sim\mathcal N(0,1)$. In the proportional asymptotic regime, if the training and test error defined in Eq.~\eqref{eq:errdef} and associated with the $M$-estimator of the problem in Eq.~\eqref{eq:ERMgen} have finite limit, they converge as in Eq.~\eqref{eq:limiterr}, with
\begin{equation}\comprimi
\begin{split}
\varepsilon_t(\alpha,\lambda,\beta_{\star}^{2})&=\textstyle\sum_c p_c\int_{\mathcal Y}\dd y\,\mathbb E_{\sigma_c,\zeta}\left[Z_0\left(y,
\bmu_c^\intercal\bbeta_\star+\frac{\sigma_cm}{\sqrt{q}}\zeta,\sigma_c^2\beta_\star^2-\frac{\sigma_c^2m^2}{q}\right)\varphi\left(y,\sigma_c\sqrt{q}\zeta+v\sigma_c^2 f_c\right)\right],\\
\varepsilon_g(\alpha,\lambda,\beta_{\star}^{2})&=\textstyle\sum_c p_c\int_{\mathcal Y}\dd y\int\dd\eta\int\dd\tau\,P_0\left(y|\tau\right)\mathbb E_{\sigma_c}\left[\mathcal N\left(\begin{psmallmatrix}{\tau}\\{\eta}\end{psmallmatrix};\begin{psmallmatrix}{
\bmu_c^\intercal\bbeta_\star}\\{t_c}\end{psmallmatrix},\sigma_c^2\begin{psmallmatrix}\beta_\star^2&m\\m&q\end{psmallmatrix}\right)\right]\varphi(y,\eta).
\end{split}\end{equation}
In particular, the estimation error is given by
\begin{equation}
\label{eq:def:esterr_k}\textstyle
\varepsilon_{\rm est}\coloneqq\lim_{d\to+\infty}\frac{1}{d}\mathbb E_{\mathcal D}\big[\|\hat\bbeta_\lambda-\bbeta_\star\|_2^2\big]=\beta_\star^2-2m+q
\end{equation}
The required set of order parameters and proximals can be found by self-consistently solving the equations below
\begin{subequations}\label{eq:sp}
\begin{equation}\comprimi
\begin{split}
	v &=\textstyle\frac{1}{d\sqrt{\hat q}}\mathbb{E}_{\bxi}[\bgg^\intercal\bxi],\\
	q &\textstyle=\frac{1}{d}\textstyle\mathbb{E}_{\bxi}[\|\bgg\|_2^2],\\
	m &=\textstyle\frac{1}{d}\mathbb{E}_{\bxi}\left[\bgg^\intercal\bbeta_\star\right],\\
	t_c &=\textstyle\mathbb{E}_{\bxi}\left[\bgg^\intercal\bmu_c\right],
\end{split}\quad
\begin{split}
\hat q_c&=\textstyle\alpha p_c \int_{\mathcal Y}\dd y\,\mathbb E_{\sigma_c,\zeta}\left[\sigma_c^2 Z_0\left(y,\bmu_c^\intercal\bbeta_\star+\frac{\sigma_cm}{\sqrt q}\zeta,\sigma_c^2\beta_\star^2-\frac{\sigma_c^2m^2}{q}\right)f_c^2\right],\\
\hat v&=\textstyle-\alpha\sum_c p_c  \int_{\mathcal Y}\dd y\,\mathbb E_{\sigma_c,\zeta}\left[\sigma_c^2 Z_0\left(y,\bmu_c^\intercal\bbeta_\star+\frac{\sigma_cm}{\sqrt q}\zeta,\sigma_c^2\beta_\star^2-\frac{\sigma_c^2m^2}{q}\right)\partial_\omega f_c\right],\\
\hat m_c&=\textstyle\alpha p_c \int_{\mathcal Y}\dd y\,\mathbb E_{\sigma_c,\zeta}\left[\sigma_c^2\partial_\mu Z_0\left(y,\bmu_c^\intercal\bbeta_\star+\frac{\sigma_cm}{\sqrt q}\zeta,\sigma_c^2\beta_\star^2-\frac{\sigma_c^2m^2}{q}\right)f_c\right],\\
\hat t_c&=\textstyle\alpha  p_c\int_{\mathcal Y}\dd y\,\mathbb E_{\sigma_c,\zeta}\left[Z_0\left(y,\bmu_c^\intercal\bbeta_\star+\frac{\sigma_cm}{\sqrt q}\zeta,\sigma_c^2\beta_\star^2-\frac{\sigma_c^2m^2}{q}\right)f_c\right],
\end{split}   
\end{equation}
where, as before, $Z_0(y,\mu,V)\coloneqq\mathbb E_z[P_0(y|\mu+\sqrt Vz)]$ and
\begin{equation}
\begin{split}
\bgg&\textstyle\coloneqq\arg\min_{\bbeta}\left(\frac{\hat v\|\bbeta\|_2^2}{2}- \sum_c\bbeta^\intercal(\hat m_c\bbeta_\star+d\hat t_c\bmu_c)-\sum_c\sqrt{\hat q_c}\bxi^\intercal\bbeta+\lambda r(\bbeta)\right)\in\R^{d},\\
f_c&\textstyle\coloneqq\arg\min_{u}\left[\frac{vu^2\sigma_c^2}{2}+\ell\left(y,t_c+\sigma_c\sqrt{q}\zeta+\sigma_c^2 v u\right)\right]\in\R.
\end{split}
\end{equation}
\end{subequations}
\end{result}

The previous set of equations covers a wide range of distributions and possibly any power-law tail decay. The generality of this setting relies on the fact that any distribution can be approximated by a possibly uncountable superposition of Gaussians \citep{Nestodoris,Alspach,ghosh2006bayesian}. Moreover, in Appendices \ref{app:universality} the following universality result, generalising results of \citet{gerace2023gaussian} and \citet{pesce2023} for the Gaussian setting, is given.
\begin{result}[Universality for uncorrelated target] Assume that, for all $c\in[K]$, $\lim_{d\to+\infty}\bmu_c^\intercal\bbeta_\star=0$ and that the labels are generated according to the model in Eq.~\eqref{eq:def:data}, the noise $\eta$ having even distribution and zero mean. Then, if the loss $\ell$ is even, the asymptotic training and generalisation errors are the same as if the covariates had distribution $p(\bx)=\mathbb E[\mathcal N(\bx;\mathbf 0,\sfrac{\sigma}{d}\bI_d)]$ with $\sigma\sim\sum_{c=1}^Kp_c\varrho_c$. Furthermore, if $\ell(y,t)=\frac{1}{2}(y-t)^{2}$ and $r(\bx)=\frac{1}{2}\|\bx\|_2^2$, under the assumption that $\mathbb E[\eta^2]<+\infty$, then the training loss is
\begin{equation}\textstyle
\lim_{\lambda\to 0^+}\frac{1}{2n}\sum_{i=1}^n(y_i-\bbeta_\lambda^\intercal\bx_i)^2\xrightarrow[\alpha=\sfrac{n}{d}=\Theta(1)]{d\to+\infty}\frac{\mathbb E[\eta^2]}{2}\left(1-\frac{1}{\alpha}\right)_+,\qquad (x)_+\coloneqq x\theta(x),    
\end{equation}
for any distribution of the random variables $\sigma_c$.
\end{result}

%% file: sections/acknowledgement.tex
The authors would like to thank Pierpaolo Vivo, Pragya Sur, Jeremias Knoblauch, Gabriel Peyr\'e, Matteo Vilucchio and Vittorio Erba for stimulating discussions and feedback. BL acknowledges funding from the \textit{Choose France - CNRS AI Rising Talents} program.

%% file: sections/appendix/replicas.tex
In this Appendix, we derive the fixed point equations for the order parameters following the analysis by \citet{loureiro2021,pesce2023,adomaityte2023} in the most general setting discussed in Result \ref{res:general}: the case in Result \ref{res:K1} is obtained by fixing $K=1$ and $\bmu_1\equiv \textbf{0}$ below, and by assuming a ridge regularisation. The dataset $\mathcal D\coloneqq \{(\bx_i,y_i)\}_{i\in[n]}$ consists of $n$ independent datapoints $\bx_i\in\R^d$ each associated to a label $y_i\in\mathcal Y$. The elements of the dataset are independently generated by using a law $P(\bx,y)$ which we assume can be put in the form of an arbitrary mixture of $K$ elliptical distributions, that we call \textit{elliptical mixture model} (EMM), parametrising the distribution of $K$ clusters $\mathcal C=\{1,\dots,K\}$,
\begin{equation}\label{app:def:datalaw}
P(\bx,y)\equiv P_0(y|\bbeta_\star^\intercal \bx)P(\bx),\qquad P(\bx)\coloneqq\sum_{c\in\mathcal C} p_c\mathbb E_{\sigma_c}\left[\mathcal N\left(\bx;\bmu_c,\sfrac{\sigma_c^2}{d}\bI_d\right)\right],  
\end{equation}
where $P_0(\bullet|\tau)$ is the distribution of the scalar label $y$ produced via the ``teacher'' $\bbeta_\star\in\R^d$ which si assumed to be sampled, once and for all, by a Gaussian distribution, $\bbeta_{\star}\sim\mathcal{N}(0,\beta_{\star}^{2}\bI_{d})$, for some $\beta_\star\in\R^*_+$. 
In the equation above, $\forall c\in\mathcal C$, $ p_c\in[0,1]$ and $\bmu_c\in\R^d$ with $\|\bmu_c\|_2^2=\Theta(\sfrac{1}{d})$. It is assumed that $\sum_c p_c=1$. For every $c$, the expectation is intended over $\sigma_c$, a positive random variable with density $\varrho_c$. We will perform our regression task searching for a set of \textit{weights} $\hat\bbeta_\lambda\in\R^d$, that will allow us to construct an estimator via a certain classifier $f\colon\R\to\mathcal Y$: 
\begin{equation}
\bx\mapsto f(\hat\bbeta_\lambda^\intercal\bx)=y,
\end{equation}
which will provide us with a prediction for a datapoint $\bx\in\R^d$. The weights will be chosen by minimising an empirical risk function in the form $\mathcal R(\bbeta)$, so that
\begin{equation}
\hat\bbeta_\lambda\coloneqq \arg\min_{\bbeta\in \R^{d}}\mathcal R(\bbeta),\qquad\text{where}\quad\mathcal R(\bbeta)\equiv \sum_{\nu=1}^n\ell\left(y_i,\bbeta^\intercal\bx_i\right)+\lambda r(\bbeta).
\end{equation}
We assume that $\ell$ is a convex loss function with respect to its second argument and $r$ is a strictly convex regularisation function: the parameter $\lambda\geq0$ will tune the strength of the regularisation. Note that this setting is slightly more general than the one given in the main text. The starting point is to reformulate the problem as an optimisation problem by introducing a Gibbs measure over the parameters $\bbeta$ depending on a positive parameter $\upbeta$,
\begin{equation}
\mu_{\upbeta}(\bbeta)\propto \e^{-\upbeta\mathcal R(\bbeta)}=
\underbrace{\e^{-\upbeta r(\bbeta)}}_{P_w}\prod\limits_{i=1}^{n}\underbrace{\exp\left[-\upbeta\ell\left(y_i,\bbeta^\intercal\bx_i\right)\right]}_{P_{\ell }},
\end{equation}
so that, in the $\upbeta\to+\infty$ limit, the Gibbs measure concentrates on $\hat\bbeta_\lambda$. The functions $P_{y}$ and $P_{w}$ can be interpreted as (unnormalised) likelihood and prior distribution respectively. Our analysis will go through the computation of the average \textit{free energy density} associated with such Gibbs measure in a specific proportional limit, i.e.,
\begin{equation}
f_\upbeta\coloneqq-\lim_{\substack{n,d\to+\infty\\\sfrac{n}{d}=\alpha}}\mathbb E_{\mathcal D}\left[\frac{\ln \mathcal Z_\upbeta}{d\upbeta}\right]=\lim_{\substack{n,d\to+\infty\\\sfrac{n}{d}=\alpha}}\lim_{s\to 0}\frac{1-\mathbb E_{\mathcal D}[\mathcal Z_\upbeta^s]}{sd\upbeta},
\end{equation}
where $\mathbb E_{\mathcal D}[\bullet]$ is the average over the training dataset, $\alpha=\Theta(1)$, and we have introduced the partition function
\begin{equation}
\mathcal Z_\upbeta\coloneqq \int \e^{-\upbeta\mathcal R(\bbeta)}\dd\bbeta.
\end{equation}

\subsection{Replica approach} 
As a first step of our replica analysis, we aim to evaluate
\begin{equation}
\mathbb E_{\mathcal D}[\mathcal Z_\upbeta^s]=\prod_{a=1}^s\int \dd \bbeta^a P_w(\bbeta^a)\left(\mathbb E_{(\bx,y)}\left[\prod_{a=1}^s P_\ell(y|\bx^\intercal\bbeta^a{})\right]\right)^n,
\end{equation}
which is the partition function replicated $s$ times and averaged over the data distribution. In the expression above we used the fact that the datapoints in our database are independently sampled. Observing that, conditioning on $c$ and $\sigma_c$, $\bx$ has Gaussian distribution, let us take the inner average introducing a new set of \textit{local fields} $\eta^a=\bx^\intercal\bbeta^a$ and $\tau=\bx^\intercal\bbeta_\star$, which are, conditional 
\begin{multline}
\mathbb E_{(\bx,y)}\left[\prod_{a=1}^s P_\ell(y|\bx^\intercal\bbeta^a)\right]
\comprimi=\mathbb E_c\left[\int_{\mathcal Y}\dd y\int_{\R^d}\dd \bx\,P_0(y|\bx^\intercal\bbeta_\star)\mathcal N(\bx;\bmu_c,\sfrac{\sigma_c^2}{d}\bI_d)\prod_{a=1}^s P_\ell(y|\bx^\intercal\bbeta^a)\right]\\\comprimi
=\mathbb E_c\left[\int\dd\bEta\int\dd\tau\int_{\mathcal Y}\dd y\,P_0\left(y|\tau\right)\prod_{a=1}^s P_\ell\left(y|\eta^a\right)\,\mathcal N\left(\begin{psmallmatrix}{\tau}\\{\bEta}\end{psmallmatrix};\begin{psmallmatrix}{\bmu_c^\intercal\bbeta_\star}\\{\bmu_c^\intercal\bbeta^a}\end{psmallmatrix},\frac{\sigma_c^2}{d}\begin{psmallmatrix}d\beta_\star^2&\bbeta_\star^\intercal\bbeta^b\\\bbeta_\star^\intercal\bbeta^a&\bbeta^a{}^\intercal\bbeta^b\end{psmallmatrix}\right)\right],
\end{multline}
where the fields are multivariate Gaussian because of the aforementioned conditional Gaussianity, and the average $\mathbb E_c[\bullet]\coloneqq\sum_c p_c\mathbb E_{\sigma_c}[\bullet]$ over the covariates' variance $\sigma_c$ is being carried over to the next steps.
We can then write 
\begin{multline}
\comprimi\mathbb E_{\mathcal D}[\mathcal Z_\upbeta^s]=\\\comprimi=\prod_{a=1}^s\int \dd \bbeta^a P_w(\bbeta^a)\left(\mathbb E_{c}\left[\int\dd\bEta\int\dd\tau\int_{\mathcal Y}\dd y\,P_0\left(y|\tau\right)\prod_{a=1}^s P_\ell\left(y|\eta^a\right)\,\mathcal N\left(\begin{psmallmatrix}{\tau}\\{\bEta}\end{psmallmatrix};\begin{psmallmatrix}{\bmu_c^\intercal\bbeta_\star}\\{\bmu_c^\intercal\bbeta^a}\end{psmallmatrix},\frac{\sigma_c^2}{d}\begin{psmallmatrix}d\beta_\star^2&\bbeta_\star^\intercal\bbeta^b\\\bbeta_\star^\intercal\bbeta^a&\bbeta^a{}^\intercal\bbeta^b\end{psmallmatrix}\right)\right]\right)^n
\\
=\prod_c\left(\prod_{a\leq b}\iint\mathcal D{\bQ}^{ab}\mathcal D\hat{\bQ}^{ab}\right)
\left(\prod_{a}\int\mathcal D{{\boldsymbol M}}^{a}\mathcal D{\hat {\boldsymbol M}}^{a}\right)\left(\prod_{a}\int\dd{t}^{a}\dd{\hat t}^{a}\right)\e^{-d\upbeta\Phi^{(s)}}.
\end{multline}
In the equation above we introduced the following \textit{order parameters}
\begin{align}
Q^{ab}_c&=\frac{\sigma_c^2}{d}\bbeta^{a\intercal}\bbeta^b\in\R,\quad a,b=1,\dots, s,\\
M^{a}_c&=\frac{\sigma_c^2}{d}\bbeta_\star^\intercal\bbeta^a\in \R,\quad a=1,\dots,s,\\
t^{a}_c&=\bmu_c^\intercal\bbeta^a\in \R,\quad a=1,\dots,s,
\end{align}
whilst the integration is over all possible order parameters, $Q^{ab}_c$ and $m^{a}_c$ to be intended as functions of $\sigma_c$. In the equation, we have also denoted the replicated free-energy
\begin{multline}\label{app:eq:replfe}
\upbeta \Phi^{(s)}(\bQ,\bhM,\bt,\bhQ,\bhM,\bht)\coloneqq\sum_c\sum_{a\leq b}\mathbb E_{\sigma_c}[\hQ^{ab}_c{Q}^{ab}_c]+\sum_c\sum_{a}\mathbb E_{\sigma_c}[\hat {M}_c^{a} M_c^{a}]+\frac{1}{d}\sum_{c,a}\hat t^a_ct^a_c\\\comprimi
-\frac{1}{d}\ln\prod_{a=1}^s\int P_w(\bbeta^a)\dd\bbeta^a\prod_c\exp\left(\sum_{a\leq b}\mathbb E_{\sigma_c}[\sigma_c^2\hQ^{ab}_c] \bbeta^{a\intercal}\bbeta^b+\sum_{a}\mathbb E_{\sigma_c}[\sigma_c^2\hat{M}^{a}_c]\bbeta^{a\intercal}\bbeta_\star+\sum_a \hat t^a_c\bbeta^{a\intercal}\bmu_c\right)\\ -\alpha\ln\mathbb E_c\left[\int\dd\bEta\int\dd\tau\int_{\mathcal Y}\dd y\,P_0\left(y|\tau\right)\prod_{a=1}^s P_\ell\left(y|\eta^a\right)\,\mathcal N\left(\begin{psmallmatrix}{\tau}\\{\bEta}\end{psmallmatrix};\begin{psmallmatrix}{t^0_c}\\{t^a_c}\end{psmallmatrix}\begin{psmallmatrix}\sigma_c^2\beta_\star^2&M^b_c\\M^a_c&Q_c^{ab}\end{psmallmatrix}\right)\right],
\end{multline}
where, for the sake of brevity, $t^0_c\coloneqq\bmu_c^\intercal\bbeta_\star$. At this point, the free energy $f_\beta$ should be computed functionally extremizing with respect to all the order parameters by virtue of the Laplace (saddle-point) approximation,
\begin{equation}
f_\upbeta=\lim_{s\to 0}\Extr_{\substack{\bQ,\bM,\bt\\\bhQ,\bhM,\bht}}\frac{\Phi^{(s)}(\bQ,\bhM,\bt,\bhQ,\bhM,\bht)}{s}.
\end{equation}
\paragraph{\bf Replica symmetric ansatz.} 
Before extremising $\Phi^{(s)}$, we make the replica symmetric assumptions, restricting therefore ourselves to a specific form of order parameters,
\begin{equation}
\begin{split}
Q^{ab}_c&=\begin{cases}R_c&a=b,\\ Q_c&a\neq b,\end{cases}\\
M_c^{a}&=M_c,\\
t_c^{a}&=t_c,
\end{split}\qquad
\begin{split}
\hQ^{ab}_c&=\begin{cases}-\frac{1}{2}\hR_c,&a=b,\\ \hQ_c&a\neq b,\end{cases}\\
\hat{M}_c^{a}&=\hat{M}_c\quad \forall a,\\
\hat{t}_c^{a}&=\hat{t}_c\quad \forall a.
\end{split}
\end{equation}
By means of this ansatz, the expressions for the fourth and fifth terms of the replicated free energy in Eq.~\eqref{app:eq:replfe} strongly simplify. The assumption is motivated by the convexity of the problem. If we denote $V_c\coloneqq R_c-Q_c$, after some work we obtain an expression for the fifth (last) term of the replicated free energy in Eq.~\eqref{app:eq:replfe}:
\begin{multline}
\ln\mathbb E_{c}\left[\int\dd\bEta\int\dd\tau\int_{\mathcal Y}\dd y\,P_0\left(y|\tau\right)\prod_{a=1}^s P_\ell\left(y|\eta^a\right)\,\mathcal N\left(\begin{psmallmatrix}{\tau}\\{\bEta}\end{psmallmatrix};\begin{psmallmatrix}{t^0_c}\\{t_c\bUno_s}\end{psmallmatrix},\begin{psmallmatrix}\sigma_c^2\beta_\star^2&M_c\bUno_s^\intercal\\M_c\bUno_s&Q_c\bUno_{s\times s}\end{psmallmatrix}\right)\right]\\\comprimi
=s\mathbb E_{c,\zeta}\left[\int_{\mathcal Y}\dd yZ_0\left(y,t^0_c+\frac{M_c\zeta}{\sqrt{Q_c}},\sigma_c^2\beta_\star^2-\frac{M_c^2}{Q_c}\right)\ln Z_\ell\left(y,t_c+\sqrt{Q_c}\zeta,V_c\right)\right]+o(s),
\end{multline}
with $\zeta\sim\mathcal N(0,1)$ and introducing the function
\begin{equation}
Z_\bullet(y,\mu,V)\coloneqq \int\frac{\dd \tau P_\bullet(y|\tau)}{\sqrt{2\pi V}}\e^{-\frac{(\tau-\mu)^2}{2V}},\qquad \bullet\in\{0,\ell\}.
\end{equation}
On the other hand, by focusing on the fourth term of the replicated free energy in Eq.~\eqref{app:eq:replfe}, writing $\hV_c=\hR_c+\hQ_c$, and introducing $\hat q_c\coloneqq \mathbb E_{\sigma_c}[\sigma_c^2\hat Q_c]$, $\hat v_c\coloneqq \mathbb E_{\sigma_c}[\sigma_c^2\hat V_c]$, and $\hat m_c\coloneqq \mathbb E_{\sigma_c}[\sigma_c^2 \hat M_c]$ we obtain
\begin{multline}\comprimi
    \frac{1}{d}\ln\prod_{a=1}^s\left(\int P_w(\bbeta^a)\dd \bbeta^a \prod_c\e^{-\frac{\hat v_c }{2}\|\bbeta^a\|_2^2+\bbeta^{a\intercal}(\hat m_c\bbeta_\star+\hat t_c\bmu_c)}\prod_{b,c}\e^{\frac{1}{2}\hat q_c\bbeta^{a\intercal}\bbeta^b}\right)=\\
    \comprimi
    =\frac{s}{d}\mathbb E_{\bxi}\ln\left[\int P_w(\bbeta)\dd\bbeta \prod_c\exp\left(-\frac{\hat v_c\|\bbeta\|_2^2}{2}+\bbeta^\intercal(\hat m_c \bbeta_\star+\hat t_c\bmu_c)+\sqrt{\hat q_c}\bxi^\intercal\bbeta\right)\right] +o(s),  
\end{multline}
where we have also introduced $\bxi\sim\mathcal N(\textbf{0},\bI_d)$. Therefore, the (replicated) \textit{replica symmetric} free-energy is given by
\begin{equation}\comprimi
\lim_{s\to 0}\frac{\upbeta}{s}\Phi^{(s)}_{\rm RS}=\frac{1}{d}\sum_{c}\hat t_ct_c+\sum_{c}\mathbb E_{\sigma_c}[\hat{M}_c M_c]+\frac{\sum_c\mathbb E_{\sigma_c}\left[\hV_c Q_c-\hQ_c V_c-\hV_c V_c\right]}{2}-\alpha\upbeta\Psi_\ell(M,Q,V)-\upbeta\Psi_w(\hat m,\hat q,\hat v)
\end{equation}
where we have defined two contributions
\begin{equation}\comprimi\begin{split}
\Psi_{\ell}(M,Q,V)&\coloneqq \frac{1}{\upbeta}\sum_c p_c\mathbb E_{\sigma_c,\zeta}\left[\int_{\mathcal Y}\dd yZ_0\left(y,t^0_c+\frac{M_c\zeta}{\sqrt{Q_c}},\sigma_c^2\beta_\star^2-\frac{M_c^2}{Q_c}\right)\ln Z_\ell\left(y,t_c+\sqrt{Q_c}\zeta,V_c\right)\right],\\
\Psi_w(\hat m,\hat q,\hat v)&\coloneqq \frac{1}{\upbeta d}\mathbb E_{\bxi}\ln\left[\int P_w(\bbeta)\dd\bbeta \prod_c\exp\left(-\frac{\hat v_c \|\bbeta\|_2^2}{2}+\bbeta^\intercal\left(\hat m_c\bbeta_\star+\hat t_c\bmu_c\right)+\sqrt{\hat q_c}\bxi^\intercal\bbeta\right)\right].
\end{split}\end{equation}
Note that we have separated the contribution coming from the chosen loss (the so-called \textit{channel} part $\Psi_\ell$) from the contribution depending on the regularisation (the \textit{prior} part $\Psi_w$). To write down the saddle-point equations in the $\upbeta\to+\infty$ limit, let us first rescale our order parameters as $\hat{M}_c\mapsto \upbeta\hat{M}_c$, $\hat t_c\mapsto d\upbeta\hat t_c$, $\hQ_c\mapsto \upbeta^2\hQ_c$, $\hV_c\mapsto \upbeta\hV_c$ and $V_c\mapsto\upbeta^{-1}V_c$. Also, for future convenience, let us rescale $Q_c\mapsto\sigma_c^2q_c$, $M_c\mapsto \sigma_c^2m_c$, $V_c\mapsto \sigma_c^2 v_c$. For $\upbeta\to+\infty$ the channel part is
\begin{equation}
\Psi_\ell(m,q,v,t)\comprimi=-\mathbb E_{c,\zeta}\left[\int_{\mathcal Y}\dd y\,Z_0\left(y,t^0_c+\frac{\sigma_cm_c\zeta}{\sqrt{q_c}},\sigma_c^2\beta_\star^2-\frac{\sigma_c^2m_c^2}{q_c}\right)\min_h\left[\frac{(h-\omega_c)^2}{2\sigma_c^2v_c}+\ell(y,h)\right]\right].
\end{equation}
where we have written $\Psi_\ell$ of a Moreau envelope and introduced for brevity
\[\omega_c\coloneqq t_c+\sigma_c\sqrt{q_c}\zeta\]
The so-called \textit{proximal} \cite{ParikhBoyd_proximalbook} will be a convenient quantity in the following
\begin{equation}
h_c\coloneqq \arg\min_{u}\left[\frac{(u-\omega_c)^2}{2\sigma_c^2v_c}+\ell(y,u)\right].
\label{app:eq:prox}
\end{equation}
A similar expression can be obtained for $\Psi_w$ by introducing the proximal
\begin{equation}\label{app:eq:general_prox_g}
\bgg=\arg\min_{\bbeta}\left(\frac{\|\bbeta\|_2^2\sum_c\hat v_c}{2} -\bbeta^\intercal\sum_c\left(\hat m_c\bbeta_\star+d\hat t_c\bmu_c\right)-\bxi^\intercal\bbeta\sum_c\sqrt{\hat q_c}+\lambda r(\bbeta)\right)\in\R^{d},
\end{equation}
we rewrite the prior part as
\begin{equation}
\Psi_w(\hat m,\hat q,\hat v,\hat t)
=-\frac{1}{d}\mathbb E_{\bxi}\left[\frac{\|\bgg\|_2^2}{2}\sum_c\hat v_c-\bgg^\intercal\sum_c\left(\hat m_c\bbeta_\star+\hat t_c\bmu_c\right)-\bxi^\intercal\bgg\sum_c\sqrt{\hat q_c}+\lambda r(\bgg)\right].
\end{equation}
The parallelism between the two contributions is evident, aside from the different dimensionality of the involved objects. The replica symmetric free energy in the $\upbeta\to+\infty$ limit is computed by extremising with respect to the introduced order parameters, 
\begin{multline}\label{eq:app:rsfree}
\comprimi f_{\rm RS}=\Extr\left[\sum_c\mathbb E_{\sigma_c}[\sigma_c^2\hat{M}_c m_c]+\frac{1}{2}\sum_c\mathbb E_{\sigma_c}\left[\sigma_c^2\left(\hV_c q_c-\hQ_c v_c\right)\right]+\sum_c t_c\hat t_c-\alpha\Psi_\ell(m,q,v,t)-\Psi_w(\hat{m},\hat q,\hat v,\hat t)\right].
\end{multline}
To do so, we have to write down a set of saddle-point equations and solve them.

\paragraph{\bf Saddle-point equations.} The saddle-point equations are derived straightforwardly from the obtained replica-symmetric free energy $f_{\rm RS}$ by functionally extremising with respect to all parameters. It is easily seen that $v_c$, $q_c$ and $m_c$ are independent from $\sigma_c$, and that it is possible to reduce the number of variables by introducing $\hat v=\sum_c\hat v_c$, so that we can write
\begin{subequations}\label{app:eq:sp}
\begin{align}
	v_c &=\frac{\mathbb{E}_{\bxi}[\bgg^\intercal\bxi]}{d\sqrt{\hat q_c}},\\
	q &=\frac{\mathbb{E}_{\bxi}[\|\bgg\|_2^2]}{d}, \\
	m &=\frac{\mathbb{E}_{\bxi}\left[\bgg^\intercal\bbeta_\star\right]}{d},\\
 t_c &=\mathbb{E}_{\bxi}\left[\bgg^\intercal\bmu_c\right].
\end{align}
\end{subequations}
which come from extremising the prior part, and the remaining equations from extremising the channel part can be rewritten as
\begin{subequations}\label{app:eq:sphat}\comprimi
\begin{align}
\hat q_c&=\alpha  p_c \int_{\mathcal Y}\dd y\,\mathbb E_{\sigma_c,\zeta}\left[\sigma_c^2 Z_0\left(y,t^0_c+\frac{\sigma_cm}{\sqrt{q}}\zeta,\sigma_c^2\beta_\star^2-\frac{\sigma_c^2m^2}{q}\right)f_c^2\right],\\
\hat v&=-\alpha \int_{\mathcal Y}\dd y\,\mathbb E_{c,\zeta}\left[\sigma_c^2 \,Z_0\left(y,t^0_c+\frac{\sigma_cm}{\sqrt{q}}\zeta,\sigma_c^2\beta_\star^2-\frac{\sigma_c^2m^2}{q}\right)\partial_\omega f_c\right],\\
\hat m_c&=\alpha  p_c\int_{\mathcal Y}\dd y\,\mathbb E_{\sigma_c,\zeta}\left[\sigma_c^2 \,\partial_\mu Z_0\left(y,t^0_c+\frac{\sigma_cm}{\sqrt{q}}\zeta,\sigma_c^2\beta_\star^2-\frac{\sigma_c^2m^2}{q}\right)f_c\right],\\
\hat t_c&=\alpha  p_c\int_{\mathcal Y}\dd y\,\mathbb E_{\sigma_c,\zeta}\left[Z_0\left(y,t^0_c+\frac{\sigma_cm}{\sqrt{q}}\zeta,\sigma_c^2\beta_\star^2-\frac{\sigma_c^2m^2}{q}\right)f_c\right],
\end{align}
\end{subequations}
with proximals defined as
\begin{equation}\begin{split}
\bgg&=\arg\min_{\bbeta}\left(\frac{\|\bbeta\|_2^2\hat v}{2} -\bbeta^\intercal\sum_c\left(\hat m_c\bbeta_\star+d\hat t_c\bmu_c\right)-\bxi^\intercal\bbeta\sum_c\sqrt{\hat q_c}+\lambda r(\bbeta)\right), \\
f_c&\coloneqq \arg\min_u\left[\frac{\sigma_c^2  v_cu^2}{2}+\ell(y,\omega_c+\sigma_c^2 v_c u)\right],\qquad \text{with}\quad\omega_c=t_c+\sigma_c\sqrt{q}\zeta.
\end{split}
\end{equation}
To obtain the replica symmetric free energy, therefore, the given set of equations has to be solved, and the result is then plugged in Eq.~\eqref{eq:app:rsfree}. The obtained saddle-point equations correspond to the ones given in Result~\ref{res:general}.

\paragraph{\bf Training and test errors.}\label{app:replica:errors}
As anticipated, the values of the order parameters, when finite, will allow to compute a number of quantities of interest. Let us start from estimating
\begin{equation}\label{app:eq:trainloss}
\varepsilon_\ell\coloneqq \lim_{n\to+\infty}\frac{1}{n}\sum\limits_{i=1}^{n}\ell(y_i, \hat\bbeta_\lambda^\intercal\bx_i).
\end{equation}
The best way to proceed is to observe that 
\begin{equation}\comprimi
\varepsilon_\ell=-\lim_{\upbeta\to+\infty}\partial_\upbeta(\upbeta\Psi_\ell)=\int_{\mathcal Y}\dd y\,\mathbb E_{c,\zeta}\left[Z_0\left(y,t^0_c+\frac{\sigma_cm}{\sqrt{q}}\zeta,\sigma_c^2\beta_\star^2-\frac{\sigma_c^2m^2}{q}\right)\ell(y,h_c)\right].
\end{equation}
This result holds for a generic function $\varphi\colon\mathcal Y\times\R\to\R$ under the hypothesis that its mean exists, so that more generally, under Assumption~\ref{ass:regime}, that is $n,d\to\infty$ jointly while keeping their ratio $\alpha=\sfrac{n}{d}$ fixed,
\begin{equation}\comprimi
\frac{1}{n}\sum_{i=1}^n\varphi(y_i,\hat\bbeta_\lambda^\intercal\bx_i)\xrightarrow{n,d\to+\infty}\int_{\mathcal Y}\dd y\,\mathbb E_{c,\zeta}\left[Z_0\left(y,t^0_c+\frac{\sigma_cm}{\sqrt{q}}\zeta,\sigma_c^2\beta_\star^2-\frac{\sigma_c^2m^2}{q}\right)\varphi(y,h_c)\right].
\end{equation}
The expressions above hold in general, but, as anticipated, important simplifications can occur in the set of saddle-point equations Eq.~\eqref{app:eq:sphat} and Eq.~\eqref{app:eq:sp} depending on the choice of the loss $\ell$ and of the regularization function $r$. The population's expectation (e.g., in the computation of the test error) of a performance function $\varphi\colon\mathcal Y\times\R\to\R$ is given instead by
\begin{equation}
\varepsilon_g\coloneqq \lim_{n\to+\infty}\mathbb E_{(\bx,y)}\left[\varphi(y,\hat\bbeta_\lambda^\intercal\bx)\right],
\end{equation}
where the expectation has to be taken on a newly sampled datapoint $(\bx,y)\not\in\mathcal D$. This expression can be rewritten as
\begin{multline}
\lim_{n\to+\infty}\mathbb E_{(\bx,y)}\left[\varphi(y,\hat\bbeta_\lambda^\intercal\bx)\right]=\lim_{n\to+\infty}\mathbb E_{y|\bx}\left[\int\dd\eta \varphi(y,\eta)\mathbb E_{\bx}\left[\delta\left(\eta-\hat\bbeta_\lambda^\intercal\bx\right)\right]\right]\\
=\int\dd\eta\int\dd\tau\int_{\mathcal Y}\dd y\,P_0\left(y|\tau\right)\varphi(y,\eta)\mathbb E_{c}\left[\mathcal N\left(\begin{psmallmatrix}{\tau}\\{\eta}\end{psmallmatrix};\begin{psmallmatrix}{t^0_c}\\{t_c}\end{psmallmatrix},\sigma_c^2\begin{psmallmatrix}\beta_\star^2&m\\m&q\end{psmallmatrix}\right)\right]
\end{multline}
whenever the quantity above is finite. This can be easily computed numerically once the order parameters are given. Finally, another relevant quantity for our investigations is the estimation mean-square error, which is finite whenever the order parameter are finite, and equal to
\begin{equation}
\varepsilon_{\rm est}\coloneqq\lim_{d\to+\infty}\frac{1}{d}\mathbb E_{\mathcal D}\big[\|\hat\bbeta_\lambda-\bbeta_\star\|_2^2\big]=\beta_\star^2-2m+q.\label{app:eq:esterr}
\end{equation}
\subsection{Bayes-optimal estimator for $K=1$}\label{app:BO}
A derivation similar to the one discussed above can be repeated to obtain information on the statistical properties of the Bayes optimal estimator presented in Result \ref{res:BO}. Given a dataset $\mathcal D$ of observations, we have that
\begin{equation}
P(\bbeta|\mathcal D)=\frac{P(\bbeta) P(\mathcal D|\bbeta)}{\mathcal Z(\mathcal D)}=\frac{P(\bbeta)}{\mathcal Z(\mathcal D)}\prod_{i=1}^nP_0(y_i|\bbeta^\intercal\bx_i)
\end{equation}
where $P(\bbeta)$ is the prior on the teacher that we assume to be Gaussian, $P(\bbeta)=\mathcal N(\bbeta;\mathbf 0,\beta_\star^2\bI_d)$, and
\begin{equation}\comprimi
\mathcal Z(\mathcal D)\coloneqq\int\dd\bbeta\,P(\bbeta) \prod_{i=1}^nP_0(y_i|\bbeta^\intercal\bx_i)=\frac{1}{(2\pi)^{n/2}}\int\dd\bbeta\,\exp\left(-\frac{\|\bbeta\|_2^2}{2\beta_\star^2}+\sum_{i=1}^n\ln P_0(y_i|\bbeta^\intercal\bx_i) \right).
\end{equation}
The calculation of $\mathcal Z(\mathcal D)$ gives access in particular to the free entropy $\phi(\mathcal D)\coloneqq \lim_{n\to\infty}\frac{1}{n}\ln\mathcal Z(\mathcal D)$. The computation of $\phi$, which has an information-theoretical interpretation as mutual information, provides us the statistics of $\bbeta$ according to the true posterior $P(\bbeta|\mathcal D)$. By assuming a concentration in the large $n$ limit, the calculation is performed on $\mathbb E_{\mathcal D}[\ln\mathcal Z(\mathcal D)]$. Using the replica trick as before, the calculation can be repeated almost identically. For the sake of simplicity, we focus on the case in which only one cluster is present, centered in the origin. It is found that the statistics of a Bayes optimal estimator can be characterised therefore by an order parameter $\sq$ satisfying a self-consistent equation not different from the one presented above (we will use below a different font to stress that we are currently analysing the Bayes optimal setting) 
\begin{equation}\label{app:eq:spbo}
\shq=\alpha  \int_{\mathcal Y}\dd y\,\mathbb E_{\sigma,\zeta}\left[\sigma^2 Z_0(y,\mu,V)\left(\partial_\mu\ln Z_0(y,\mu,V)\right)^2\Big|_{\substack{\mu=\sigma\sqrt{\sq}\zeta\\ V=\sigma^2(\beta_\star^2-\sq)}}\right],\quad 
	\sq =\frac{\beta_\star^4\shq}{1+\beta_\star^2\shq}.
\end{equation}
with $Z_0(y,\mu,V)\coloneqq \mathbb E_z[P_0(y|\mu+\sqrt V z)]$ with $z\sim\mathcal N(0,1)$. We can then compute the Bayes optimal estimation error for the Bayes optimal estimator $\hat\bbeta_{\rm BO}=\mathbb E_{\bbeta|\mathcal D}[\bbeta]$ as
\begin{equation}
\varepsilon_{\rm BO}=\lim_{d\to+\infty}\frac{1}{d}\|\bbeta_\star-\hat\bbeta_{\rm BO}\|_2^2=\beta_\star^2-\sq.
\end{equation}

%% file: sections/appendix/ridge.tex
{Let us now fix $\ell_2$ regularisation $r(\bx)=\frac{1}{2}\|\bx\|_2^2$. In this case, the computation of $\Psi_w$ can be performed explicitly via a Gaussian integration, and the saddle-point equations can take a more compact form that is particularly suitable for a numerical solution. In particular, the prior proximal with general expression in Eq.~\eqref{app:eq:general_prox_g} simplifies to
\begin{equation}
\bgg
= \frac{\sum_c\left(\hat m_c\bbeta_\star+d\hat t_c\bmu_c\right)+\sum_c\sqrt{\hat q_c}\bxi}{\lambda+\hat v}
\end{equation}
so that the prior (non-hat) saddle-point equations obtained from $\Psi_w$ become
\begin{equation}\label{app:eq:ridge}\comprimi
\begin{split}
	q &\textstyle=\frac{1}{d}\left(\sum_c\frac{\hat m_c\bbeta_\star+d\hat t\bmu_c}{\lambda+\hat v}\right)^2+\left(\sum_c\frac{\sqrt{q_c}}{\lambda+\hat v}\right)^2\\
	m &\textstyle= \frac{\sum_c\left(\beta_\star^2\hat m_c+t^0_c\hat t_c\right)}{\lambda+\hat v}\\
	v_c &\textstyle=\frac{1}{\lambda+\hat v}\sum_{c'}\sqrt{\frac{\hat q_{c'}}{\hat q_c}}\\
	t_c &\textstyle= \frac{\sum_{c'}\left(\hat t_{c'}\mu_{c'c}+t^0_c\hat m_{c'}\right)}{\lambda+\hat v},
\end{split}\quad 
\begin{split}
\hat q_c&\textstyle=\alpha  p_c \int_{\mathcal Y}\dd y\,\mathbb E_{\sigma_c,\zeta}\left[\sigma_c^2 Z_0\,f_c^2\right],\\
\hat v&\textstyle=-\alpha \int_{\mathcal Y}\dd y\,\mathbb E_{c,\zeta}\left[\sigma_c^2 \,Z_0\,\partial_\omega f_c\right],\\
\hat m_c&\textstyle=\alpha  p_c\int_{\mathcal Y}\dd y\,\mathbb E_{\sigma_c,\zeta}\left[\sigma_c^2 \,\partial_\mu Z_0\,f_c\right],\\
\hat t_c&\textstyle=\alpha  p_c\int_{\mathcal Y}\dd y\,\mathbb E_{\sigma_c,\zeta}\left[Z_0\,f_c\right]
\end{split}
\end{equation}
We have used the shorthand notation $Z_0\equiv Z_0\left(y,t^0_c+\frac{\sigma_cm}{\sqrt{q}}\zeta,\sigma_c^2\beta_\star^2-\frac{\sigma_c^2m^2}{q}\right)$ and $\mu_{cc'}\coloneqq d\bmu_{c'}^\intercal\bmu_c$.
\paragraph{\bf Regression on one cloud: consistency} In the special case in which $|\mathcal C|=1$ and the cloud is centered in the origin, we have $t_1=\hat t_1=0$ and, dropping the subscript referring to the cluster, the equations are
\begin{equation}\label{app:sp:ridge1}
\begin{split}
	q &=\frac{\beta_\star^2\hat m^2+\hat q}{(\lambda+\hat v)^2}\\
	m &= \frac{\beta_\star^2\hat m}{\lambda+\hat v}\\
	v &=\frac{1}{\lambda+\hat v},
\end{split}\quad
\begin{split}
\hat q&\textstyle=\alpha  \int_{\mathcal Y}\dd y\,\mathbb E_{\sigma,\zeta}\left[\sigma^2 Z_0\left(y,\frac{\sigma m}{\sqrt{q}}\zeta,\sigma^2\beta_\star^2-\frac{\sigma^2m^2}{q}\right)f^2\right],\\
\hat v&\textstyle=-\alpha \int_{\mathcal Y}\dd y\,\mathbb E_{\sigma,\zeta}\left[\sigma^2\,Z_0\left(y,\frac{\sigma m}{\sqrt{q}}\zeta,\sigma^2\beta_\star^2-\frac{\sigma^2m^2}{q}\right)\partial_\omega f\right],\\
\hat m&\textstyle=\alpha\int_{\mathcal Y}\dd y\,\mathbb E_{\sigma,\zeta}\left[\sigma^2\,\partial_\mu Z_0\left(y,\frac{\sigma m}{\sqrt{q}}\zeta,\sigma^2\beta_\star^2-\frac{\sigma^2m^2}{q}\right)f\right],
\end{split}
\end{equation}
where as usual $f\coloneqq \arg\min_u\left[\frac{\sigma^2 vu^2}{2}+\ell(y,\omega+\sigma^2 v u)\right]$ and $\omega=\sigma\sqrt{q}\zeta$. Let us now perform the rescaling $v\mapsto \alpha v$, $\hat q\mapsto \alpha \hat q$, $\hat m\mapsto \alpha \hat m$, and $\hat v\mapsto \alpha \hat v$, where $v=O(1)$, $\hat v=O(1)$, $\hat m=O(1)$, $\hat q=O(1)$. Then, under these assumptions, in the $\alpha\to+\infty$ limit
\begin{equation}\label{app:sp:ridgeres}
\begin{split}
	q &=\frac{\beta_\star^2\hat m^2}{\hat v^2}\\
	m &= \frac{\beta_\star^2\hat m}{\hat v}\\
	v &=\frac{1}{\hat v},
\end{split}\quad
\begin{split}
\hat q&\textstyle=\int_{\mathcal Y}\dd y\,\mathbb E_{\sigma,\zeta}\left[\sigma^2 Z_0\left(y,\frac{\sigma m}{\sqrt{q}}\zeta,\sigma^2\beta_\star^2-\frac{\sigma^2m^2}{q}\right)f^2\right],\\
\hat v&\textstyle=-\int_{\mathcal Y}\dd y\,\mathbb E_{\sigma,\zeta}\left[\sigma^2\,Z_0\left(y,\frac{\sigma m}{\sqrt{q}}\zeta,\sigma^2\beta_\star^2-\frac{\sigma^2m^2}{q}\right)\partial_\omega f\right],\\
\hat m&\textstyle=\int_{\mathcal Y}\dd y\,\mathbb E_{\sigma,\zeta}\left[\sigma^2\,\partial_\mu Z_0\left(y,\frac{\sigma m}{\sqrt{q}}\zeta,\sigma^2\beta_\star^2-\frac{\sigma^2m^2}{q}\right)f\right],
\end{split}
\end{equation}
Moreover, in the large $\alpha$ limit, $f=-\partial_\eta\ell(y,\eta)|_{\eta=\sigma\sqrt q\zeta}$. It follows that, independently from the adopted (convex) loss, the angle $\pi^{-1}\arccos\frac{m}{\beta_\star\sqrt q}$ between the estimator $\hat\bbeta_\lambda$ and the teacher $\bbeta_\star$ goes to zero as $\alpha\to+\infty$. In this limit, it is easily found that $\varepsilon_{\rm est}\to \beta_\star^{-2}(m-\beta_\star^2)^2$, hence the estimator is consistent if $m\to \beta_\star^2$.

\paragraph{\bf Uncorrelated teachers: universality}\label{app:universality}
To study the universality properties in the ridge setting, let us introduce two new assumptions. The first one expresses in general the fact that the teacher $\bbeta_\star$ is completely uncorrelated from the centroids $\bmu_c$ of our datasets.
\begin{assumption}\label{app:assump1}
For all $c\in[K]$, $\lim_{d\to+\infty}t_c^0=0$.
\end{assumption}
\noindent This assumption holds, for example, by assuming centroids such as $\bmu_c\sim\mathcal N(\mathbf 0,\sfrac{1}{d}\bI_d)$. We add, on top of this, a symmetry assumption on the probabilistic law of the label generation, and on the loss.

\begin{assumption}\label{app:assump2}The following symmetry properties hold
\begin{equation}
P_0(y|\tau)=P_0(-y|-\tau),\qquad \ell(y,\eta)=\ell(-y,-\eta).  
\end{equation}
\end{assumption}

Under Assumption \ref{app:assump1} and Assumption \ref{app:assump2}, $t_c=\hat t_c=0$ $\forall c$ is a saddle-point solution of the Eqs.~\eqref{app:eq:ridge}. Indeed, if $\hat t_c=0$ the prior equation implies $t_c=0$. On the other hand, if $t_c=0$, $\hat t_c=0$, due to the parity assumption. We recover therefore in our setting the \textit{mean universality} discussed by \citet{pesce2023} in the Gaussian setting: under the assumptions above, the learning task is mean-independent and equivalent to a task performed on $c$ clouds all centered in the origin, i.e., a problem obtained assuming $\bx\sim\mathbb E_{c}[\mathcal N(\mathbf 0,\sfrac{\sigma_c}{d}\bI_d)]$. Note that in the Gaussian setting --- recovered in our setting by assuming a Dirac-delta density $\varrho_c(\sigma)=\delta(\sigma-\bar\sigma_c)$ for some fixed $\bar\sigma_c$ for each $c\in[C]$ --- \citet{pesce2023} observed that \textit{covariance universality} holds for $\lambda\to 0^+$: namely, $\varepsilon_g$ and $\varepsilon_\ell$ are independent on the covariance of the clouds. This fact does not extend to the case in which the distribution of $\sigma_c$ is not atomic (not even in the case in which $\sigma_c$ are all identically distributed), as it has been verified by \citet{adomaityte2023}.

\subsection{Square loss on a linear model}\label{app:square}
If we consider a square loss $\ell(y,\eta)=\frac{1}{2}\left(y-\eta\right)^2$, then an explicit formula for the proximal can be found, namely
\begin{equation}
f_c=\frac{y-\omega_c}{1+v_c\sigma_c^2},\qquad \omega_c=t_c+\sigma_c\sqrt{q}\zeta,
\end{equation}
so that the second set (hat-parameters) of saddle-point equations  Eq.~\eqref{app:eq:sphat} can be made even more explicit. Let us assume furthermore that labels are generated according to the linear model in Eq.~\eqref{eq:def:data}, where the noise term $\eta_i$ has $\mathbb E[\eta_i]=0$ and $\hat\sigma^2_0\coloneqq\mathbb E[\eta_i^2]<+\infty$. In this setting, the channel equations can be written as
\begin{subequations}\label{app:se:ridge}
\begin{align}
\hat q_c&=\alpha p_c\hsigma^2_0 \mathbb E_{\sigma_c}\left[\frac{\sigma_c^2}{(1+v_c\sigma_c^2)^2}\right]+\alpha p_c \mathbb E_{\sigma_c}\left[\left(\frac{\sigma_c^2}{1+v_c\sigma_c^2  }\right)^2\right](\beta_\star^2-2 m+q+(t_c^0-t_c)^2),\\
\hat v&=\alpha \mathbb E_{c}\left[\frac{\sigma_c^2}{1+ v_c\sigma_c^2}\right],\\
\hat m_c&=\alpha p_c\mathbb E_{\sigma_c}\left[\frac{\sigma_c^2}{1+v_c\sigma_c^2}\right]\\
\hat t_c&=\alpha p_c(t_c^0-t_c)\mathbb E_{\sigma_c}\left[\frac{1}{1+v_c\sigma_c^2}\right].
\end{align}
\end{subequations}
Note that if $\mathbb E[\eta_i^2]$ is not finite, the expressions become meaningless. Moreover, the dependence on the noise distribution is through its second moment only: the possible power-law behavior of the noise distribution \textit{does not} affect the order parameters. On the other hand, the equations are \textit{not} affected by the possible divergence of $\mathbb E[\sigma_c^2]$, despite the fact that, in this case, the average loss is not finite. We interpret this fact as a reflection of the existence of the \textit{order parameters} associated with the minimizer. The meaningfulness of such parameters beyond the regime in which the asymptotic average loss is finite has been verified numerically. 

In this setting the generalisation error becomes
\begin{equation}
\varepsilon_g\coloneqq\lim_{d\to+\infty}\mathbb E\left[\left(y-\hat\bbeta_\lambda^\intercal\bx\right)^2\right]=\hsigma_0^2+\sum_{c}p_c(t_c^0-t_{c})^2+(\beta_\star^2-2m+q)\sum_c p_c\mathbb E_{\sigma_c}[\sigma^2_c],
\end{equation}
which, as expected, is finite if and only if the variance of the covariates is finite $\mathbb E_{\sigma_c}[\sigma^2_c]<+\infty$ for all clouds $c$. The training loss, on the other hand, is
\begin{multline}
\varepsilon_\ell\coloneqq\lim_{d\to+\infty}\frac{1}{2n}\sum_{i=1}^n\left(y_i-\hat\bbeta_\lambda^\intercal\bx_i\right)^2\xrightarrow{d\to+\infty}\\\frac{1}{2}\mathbb E_{c}\left[\frac{\hsigma_0^2+(t_c^0-t_c)^2}{(1+v_c\sigma_c^2)^2}\right]+\frac{1}{2}\mathbb E_{c}\left[\frac{\sigma_c^2}{(1+v_c\sigma_c^2 )^2}\right](\beta_\star^2-2m+q).\end{multline}

\paragraph{\bf Strong universality of $\varepsilon_\ell$ for $\lambda\to 0^+$} 
We will now show that, under Assumption \ref{app:assump1} (Assumption \ref{app:assump2} is automatically satisfied in the setting under consideration), a strong universality of the training loss holds. Let us put ourselves in the case of a single cluster centered at the origin (an assumption that is not restrictive, as shown above). In this case, let us introduce
\begin{equation}
Y(v)\coloneqq v\mathbb E_\sigma\left[\frac{\sigma^2}{1+v\sigma^2}\right]    
\end{equation}
which can be interpreted in terms of the Stieltjes transform of the random variable $\sigma^2$. The saddle-point equations can be rewritten in terms of this function as
\begin{equation}\label{app:eq:spC1}
\begin{split}
\hat q&=\alpha \hsigma_0^2\partial_vY(v)+\alpha(Y(v)-v\partial_v Y(v))\frac{\beta_\star^2-2 m+q}{v^2},\\
\hat m&=\alpha\frac{Y(v)}{v},
\end{split}\quad
\begin{split}
	q &=(\hat m^2\beta_\star^2+\hat q)v^2,\\
	m &= \beta_\star^2\hat m v,
\end{split}
\end{equation}
with an additional fixed-point equation for $v$,
\begin{equation}
1-\lambda v=\alpha Y(v).
\end{equation}
A little bit of manipulation leads to expressing the estimation error as
\begin{equation}\varepsilon_{\rm est}=v\left(\hat\sigma_0^2+\frac{(\beta^2_\star\lambda-\hat\sigma_0^2)\lambda}{\alpha\partial_vY(v)+\lambda}\right).\label{eq:app:ests}\end{equation}
Analogously, the training loss can be written as 
\begin{equation}
\varepsilon_\ell=\frac{\hsigma_0^2(1-Y(v)-v\partial_v Y(v))}{2}+\frac{(\beta_\star^2-2m+q)\partial_v Y(v)}{2}.  
\end{equation}
In the limit $\lambda\to 0$,
\begin{equation}
x\coloneqq \frac{\beta_\star^2-2m+q}{v}=\frac{(Y(v)-v\partial_v Y(v))x+\hsigma_0^2v\partial_v Y(v)}{Y(v)}\Rightarrow x=\hsigma_0^2
\end{equation}
so that $\varepsilon_\ell=\frac{1}{2}(1-Y(v))\hsigma_0^2$. The quantity $Y(v)$ can be extracted from the equation for $v$, as it has to satisfy, in the zero regularisation limit, $\alpha\,Y(v)=1\Rightarrow Y(v)=\frac{1}{\alpha}$ which is a valid solution for $\alpha>1$ only. As a result, we obtain a \textit{universal} formula for the training loss for $\alpha>1$ already identified by \citet{pesce2023} for the purely Gaussian case is found,
\begin{equation}
\varepsilon_\ell=\frac{\hsigma_0^2}{2}\left(1-\frac{1}{\alpha}\right)_+,\qquad\text{where}\quad  (x)_+\coloneqq x\theta(x).
\end{equation}
Note that the formula above is valid for \textit{any} distribution of $\sigma$, including distributions with no second moment.

\paragraph{Estimation error rates} We conclude this section by extracting the estimation error rate for $n\to+\infty$ and large but fixed $d$, i.e., for $\alpha\to+\infty$. For simplicity, let us focus, once again, on the case $K=1$ and $\bmu_1=\mathbf 0$, corresponding to the fixed-point equations given in Eq.~\eqref{app:eq:spC1}. From Eq.~\eqref{app:eq:spC1} $v$ satisfies the equation $\alpha Y(v)=1-\lambda v$: as $Y(v)\in[0,1]$ and $v>0$, then for $\alpha\to+\infty$ we must have $Y(v)\to 0$ and $v\to 0$, so that for $\alpha\to+\infty$, $Y(v)=\frac{1}{\alpha}+O(\alpha^{-1})$. In this limit, therefore, by direct inspection of the fixed-point equations, $q\to \beta_\star^2$ and $m\to \beta_\star^2$ so that $\varepsilon_{\rm est}\to 0$ and the estimator $\hat\bbeta_\lambda$ is unbiased. Moreover, due to Eq.~\eqref{eq:app:ests}, the scaling of $\varepsilon_{\rm est}$ is the scaling of $v$ with $\alpha$.

Let us now extract the asymptotic rates of the estimation error. In the following, we assume a power-law tail decay of the standard deviation of the covariates, with density parametrised by $a$ such that $\varrho(\sigma)\sim\sigma^{-2a-1}$ for $\sigma\gg 1$. We find three value ranges of the parameter $a$, each corresponding to a different scaling of the estimation error when $\alpha\to\infty$.
\begin{description}
\item[Finite covariate variance $a>1$] In the hypothesis that $\sigma_0^2\coloneqq \mathbb E[\sigma^2]$ is finite, and thus the covariates have finite variance, the expansion of $Y(v)$ for small $v$ is $Y(v)\simeq v\sigma_0^2+o(v)$, and it is found that
    \begin{equation}
    q=\beta_\star^2+\frac{\hsigma_0^2-2\beta_\star^2\lambda}{\sigma_0^2}\frac{1}{\alpha}+o\left(\frac{1}{\alpha}\right),\qquad m=\beta_\star^2-\frac{\lambda\beta_\star^2}{\sigma_0^2}\frac{1}{\alpha}+o\left(\frac{1}{\alpha}\right),\qquad v=\frac{1}{\sigma_0^2\alpha}+o\left(\frac{1}{\alpha}\right),
    \end{equation}
    which, plugged into our general formula for $\varepsilon_{\rm est}$, imply a scaling $\varepsilon_{\rm est}\sim\alpha^{-1}$ for large $\alpha$.
    \item[Infinite covariate variance $0<a<1$] If $a\in(0,1)$, then $\mathbb E[\sigma^2]=+\infty$ and $Y(v)$ has an expansion for small $v$ in the form $Y(v)=\tsigma_0^2 v^a+O(v)$ for some finite positive quantity $\tsigma_0^2$. Such asymptotic implies that $v\simeq (\tsigma_0^2\alpha)^{-\sfrac{1}{a}}$ for $\alpha\gg 1$. By replacing this in the fixed point equations, it is found that $m=\beta_\star^2-\lambda\beta_\star^2(\tsigma_0^2\alpha)^{-\sfrac{1}{a}}+o(\sfrac{1}{\alpha})$ and $q=\beta_\star^2+(\hsigma_0^2-2\lambda\beta_\star^2)(\tsigma_0^2\alpha)^{-\sfrac{1}{a}}+o(\sfrac{1}{\alpha})$, so that the estimation error scales as $\varepsilon_{\rm est}\sim \alpha^{-\sfrac{1}{a}}$.
    \item[The marginal $a=1$ case] The $a=1$ case also corresponds to an infinite covariate variance. However, in this case, as $v\to 0$, $Y(v)\simeq \tsigma_0^2 v\log v$ for some positive constant $\tsigma_0^2$. Therefore $v=(\tsigma_0^2\alpha\ln \alpha)^{-1}$. Consequently, for $\alpha\gg 1$ $m\simeq \beta_\star^2-\lambda\beta_\star^2(\tsigma_0^2\alpha\ln \alpha)^{-1}$ and $q\simeq \beta_\star^2+\left(\hsigma_0^2-2\lambda\beta_\star^2\right)(\tsigma_0^2\alpha\ln \alpha )^{-1}$, so that $\varepsilon_{\rm est}\sim (\alpha\ln\alpha)^{-1}$.
\end{description}
The discussion above has been summarized in the first part of the Result \ref{res:scaling} in the main text.


\subsection{Huber loss and robust regression}
A toy model for the study of robustness has been recently introduced by \citet{Vilucchio2023}: in this Appendix, we will adopt the setup investigated therein, extending it to the case of fat tails. We will focus on the one cloud case, $K=1$, so that $P(\bx)=\mathbb E_\sigma[\mathcal N(\bx;\mathbf 0,\sfrac{\sigma^2}{d}\bI_d)]$, and $P_0(y|\tau)=\mathbb E[\mathcal N(y;\vartheta \tau,\hsigma^2)]$, where the expectation is taken over the joint distribution $\hat\varrho$ for the pair $(\vartheta,\hsigma)$ of (possibly correlated) random variables. \citet{Vilucchio2023} adopted, in particular, a distribution $\hat\varrho(\vartheta,\hsigma)=\epsilon\delta_{\vartheta,\vartheta_{\rm out}}\delta_{\hsigma,\hsigma_{\rm out}}+(1-\epsilon)\delta_{\vartheta,1}\delta_{\hsigma,\hsigma_{\rm in}}$ for $\epsilon\in[0,1]$, with $(\vartheta_{\rm out},\hsigma_{\rm out})$ referring to ``outlier labels'', and $(1,\hsigma_{\rm in})$ referring to ``inlier labels'', and $\delta_{a,b}$ being the Dirac-delta zero-measure density that evaluates to $1$ when $a=b$ and is zero everywhere else. The general fixed-point equations derived above can be adapted to this case quite easily. 

\paragraph{\bf A different estimator}Before considering the application of our approach to this model, let us comment on the fact that, as observed elsewhere \cite{Vilucchio2023}, in this case, it can be interesting to consider, beyond the ERM estimator $\hat\bbeta_\lambda$ and the Bayes-optimal estimator $\hat\bbeta_{\rm BO}$,  the estimator that minimises the (posterior-averaged) mean-square test error
\begin{equation}
\hat\bbeta_{g,{\rm BO}}=\arg\min_{\bbeta}\mathbb E_{\hat\bbeta|\mathcal D}\left[\mathbb E_{(y,\bx)|\hat\bbeta}\left[\left(y-\bbeta^\intercal\bx\right)^2\right]\right]=\mathbb E[\vartheta]\hat\bbeta_{\rm BO}.
\end{equation}
In the expression above, $\mathbb E_{(y,\bx)|\hat\bbeta}$ expresses the fact that the pair $(y,\bx)$ has been generated with a teacher vector $\hat\bbeta$, sampled by the posterior. Using the results on the Bayes optimal estimator, it is simple to derive the errors obtained by using $\hat\bbeta_{g,{\rm BO}}$ under the assumptions of finite covariates covariance, $\sigma_0^2\coloneqq\mathbb E[\sigma^2]<+\infty$, and finite noise variance, $\hsigma_0^2\coloneqq\mathbb E[\hsigma^2]$,
\[\varepsilon_{g,{\rm BO}}\coloneqq\mathbb E_{(y,\bx)}\left[\left(y-\bx^\intercal\bbeta_{g,{\rm BO}}\right)^2\right]=\sigma^2_0\left(\beta_\star^2\mathbb E[\vartheta^2]-\mathbb E[\vartheta]^2\sq\right)+\hsigma_0^2,\]
where $\sq$ is provided by Eq.~\eqref{app:eq:spbo}. As in the pure Gaussian case, by imposing the ansatz $\sq=\beta_\star^2-\frac{1}{\alpha}q_0+\Theta(\alpha^{-2})$, and consequently $\shq=\alpha \hat q_0+\Theta(1)$ for large $\alpha$, we can obtain
\begin{equation}
\frac{1}{q_0}=\hat q_0=\int_{\mathcal Y}\dd y\,\mathbb E_{\sigma,\zeta}\left[\sigma^2 P_0(y|\omega)\left(\partial_\omega\ln P_0(y|\omega)\right)^2\Big|_{\omega=\sigma\beta_\star\zeta}\right].
\end{equation}
In the $\alpha\to +\infty$ limit, then, $\sq\to\rho$ and $\varepsilon_{\rm est,BO}\coloneqq\lim_{d\to+\infty}\frac{1}{d}\mathbb E_{\mathcal D}[\|\hat\bbeta_{g,{\rm BO}}-\bbeta\|_2^2]=\frac{q_0}{\alpha}+\Theta(\alpha^{-2})\to 0$. On the other hand, $\varepsilon_{g,{\rm BO}}=\hsigma_0^2+\sigma_0^2\beta_\star^2\mathrm{Var}[\vartheta]-\frac{\sigma^2_0 q_0}{\alpha}+\Theta(\alpha^{-2})$.

\subsubsection{Huber loss}
The Huber loss is a strongly convex loss depending on a tunable parameter $\delta\geq 0$ and is defined as
\begin{equation}\label{eq:app:huberloss}
\ell_\delta(y,\eta)=
\begin{cases}
\frac{(y-\eta)^2}{2}&\text{if $|y-\eta|<\delta$}\\
\delta|y-\eta|-\frac{\delta^2}{2}&\text{otherwise}.
\end{cases}
\end{equation}
This loss is widely adopted in robust regression as it is less sensitive to outliers than the most commonly adopted square loss. In the most general, multicluster setting discussed in Appendix \ref{app:replica}, its associated saddle-point equations can be obtained by using the following expression for the associated proximal
\begin{equation}\comprimi
h_c=\omega_c+\frac{(y-\omega_c)v_c\sigma_c^2}{\max(\delta^{-1}|y-\omega_c|,1+v_c\sigma_c^2)}\Leftrightarrow f_c=\frac{y-\omega_c}{\max(\delta^{-1}|y-\omega_c|,1+v_c\sigma_c^2)},\quad\text{where}\quad \omega_c=t_c+\sigma_c\sqrt{q}\zeta.
\end{equation}
Using the usual ridge regularisation, $r(\bx)=\frac{1}{2}\|\bx\|_2^2$, the prior equations (non-hat parameters) are the same in Eq.~\eqref{app:sp:ridge1}. The channel equations --- that depend on the choice of the loss --- are instead
\begin{subequations}\label{app:sp:huber}
\begin{align}
\hat m&=\alpha\mathbb E\left[\frac{\sigma^2\vartheta\erf  \chi}{1+v\sigma^2}\right]\\    
\hat q&=\alpha\mathbb E\left[\frac{\sigma^2\psi \erf \chi}{(1+v\sigma^2)^2}+\sigma^2\delta^2(1-\erf \chi)-\sqrt{\frac{2\psi}{\pi}}\frac{\sigma^2\delta\e^{-\chi^2}}{1+v\sigma^2}\right]\\  
\hat v&=\alpha\mathbb E\left[\frac{\sigma^2\erf \chi}{1+v\sigma^2}\right].
\end{align}    
\end{subequations}
where the expectation is over all random variables involved in the expressions (namely, $\sigma$, $\hsigma$, and $\vartheta$) and we used the short-hand notation
\begin{equation}
\psi\coloneqq \hsigma^2+\sigma^2(\vartheta^2\beta_\star^2-2\vartheta m+q),\qquad \chi\coloneqq \frac{\delta(1+v\sigma^2)}{\sqrt{2\psi}}.
\end{equation}

With the usual notation convention $\hsigma_0^2\coloneqq\mathbb E[\hsigma^2]$ and $\sigma^2_0\coloneqq\mathbb E[\sigma^2]$, the estimation error is given by the general formula in Eq.~\eqref{app:eq:esterr}, whereas the generalisation error is
\begin{equation}\varepsilon_g\coloneqq\mathbb E\left[(y-\hat\bbeta_\lambda^\intercal\bx)^2\right]=\hsigma_0^2+(\beta_\star^2\mathbb E[\vartheta^2]-2\mathbb E[\vartheta]m+q)\sigma^2_0,
\end{equation}
$\varepsilon_g$ being finite if $\sigma^2<+\infty$, $\hsigma_0^2<+\infty$ and $\mathbb E[\vartheta^2]<+\infty$. 

\paragraph{\bf Asymptotic Bayes optimal performances} We aim now at extrapolating the large-$\alpha$ behavior of such errors and at studying the consistency of $\hat\bbeta_\lambda$ with respect to the Bayes optimal estimators discussed in Section~\ref{app:BO}. To do so, we rescale $\hat m\mapsto\alpha\hat m$, $\hat v\mapsto\alpha\hat v$, $v\mapsto\alpha^{-1} v$ and $\hat q\mapsto\alpha\hat q$. We also assume that $\lambda\mapsto\lambda+\alpha\lambda'$ (the role of $\lambda'\neq 0$ will be clear in the following). The set of fixed point equations become, for $\alpha\to+\infty$
\begin{equation}
\begin{split}
\hat m&=\mathbb E\left[\sigma^2\vartheta\erf  \bar\chi\right]\\    
\hat q&=\mathbb E\left[\sigma^2\psi \erf \bar\chi+\sigma^2\delta^2(1-\erf \bar\chi)-\sqrt{\frac{2\psi}{\pi}}\sigma^2\delta\e^{-\bar\chi^2}\right]\\  
\hat v&=\mathbb E\left[\sigma^2\erf \bar\chi\right].
\end{split}
\qquad 
\begin{split}
	q &=\frac{\beta_\star^2\hat m^2}{(\lambda'+\hat v)^2}\\
	m &= \frac{\beta_\star^2\hat m}{\lambda'+\hat v}\\
	v &=\frac{1}{\lambda'+\hat v}
\end{split},\qquad \bar\chi\coloneqq\frac{\delta}{\sqrt{2\psi}}.
\end{equation}
In this limit, as $\beta_\star^2 q=m^2$, $\psi=\hsigma^2_0+\frac{\sigma^2_0}{\beta_\star^2}(m-\beta_\star^2\vartheta)^2$, we find
\begin{equation}
\varepsilon_{\rm est}=\frac{(m-\beta_\star^2)^2}{\beta_\star^2},\qquad \varepsilon_g=\hsigma_0^2+\sigma^2\frac{\mathbb E[(m-\vartheta\beta_\star^2)^2]}{\beta_\star^2}.
\end{equation}
It is possible to choose $\lambda'$ so that $\lim_{\alpha\to+\infty}\varepsilon_g=\lim_{\alpha\to+\infty}\varepsilon_g^{\rm BO}$, i.e.,
\begin{equation}
\hsigma_0^2+\sigma^2_0\frac{\mathbb E[(m-\vartheta\beta_\star^2)^2]}{\beta_\star^2}=\hsigma_0^2+\sigma^2_0\beta_\star^2\mathrm{Var}[\vartheta]\Rightarrow m=\beta_\star^2\mathbb E[\vartheta].
\end{equation}
We can try to satisfy this condition by tuning properly $\lambda'$, under the constraint that $\lambda'\geq 0$. We can write in particular
\begin{equation}
\lambda'=\frac{\hat m}{\mathbb E[\vartheta]}-\hat v=\frac{\mathbb E[\sigma^2\vartheta\erf\bar\chi]-\mathbb E[\vartheta]\mathbb E[\sigma^2\erf\bar\chi]}{\mathbb E[\vartheta]}\geq 0\Rightarrow \mathbb E[\sigma^2(\vartheta-\mathbb E[\vartheta])\erf\bar\chi]\geq 0\end{equation}
to be computed with
\begin{equation}
\bar\chi\equiv\frac{\delta}{\sqrt{2\psi_g}},\quad\psi_g=\hsigma^2_0+\sigma^2_0\beta_\star^2(\mathbb E[\vartheta]-\vartheta)^2.
\end{equation}
Note that the condition is always satisfied in the case of the square loss (i.e., for $\delta\to+\infty\Leftrightarrow\bar\chi\to 1$).

\paragraph{\bf Consistency of the estimator} The consistency of the estimator can be imposed by properly tuning $\lambda$, by requiring that $\lim_{\alpha\to +\infty}\varepsilon_{\rm est}=0$, i.e., $m=\beta_\star^2$ in this limit. In the same spirit as above, this implies a condition on $\lambda'$ given by
\begin{equation}
\lambda'=\hat m-\hat v=\mathbb E[\sigma^2\vartheta\erf\bar\chi]-\mathbb E[\sigma^2\erf\bar\chi]\geq 0\Rightarrow \mathbb E[\sigma^2(\vartheta-1)\erf\bar\chi]\geq 0
\end{equation}
to be computed with
\begin{equation}
\bar\chi\equiv\frac{\delta}{\sqrt{2\psi_{\rm est}}},\quad\psi_{\rm est}=\hsigma^2_0+\sigma^2_0\beta_\star^2(1-\vartheta)^2.
\end{equation}
When imposing the equality, the conditions above provide the values of $\delta$ (if any) for a consistent estimator if $\lambda=\Theta(1)$ in the $\alpha\to+\infty$ limit.

\paragraph{\bf Estimation error rates} Let us focus now on the estimation error rates for the Huber loss. With respect to the model above, adopt $\hat\varrho(\vartheta,\sigma)=\delta_{\vartheta,1}\varrho(\sigma)$. We will also assume the noise to be Gaussian with variance $\hat\sigma_0^2$. We can rewrite the equations in terms of $$\comprimi\Upsilon(v)\coloneqq v\mathbb E\left[\frac{\sigma^2\erf\chi}{1+v\sigma^2}\right],\quad\text{with}\quad\chi\coloneqq\frac{\delta(1+v\sigma^2)}{\sqrt{2\psi}},\quad\psi=\hat\sigma_0^2+\sigma^2\varepsilon_{\rm est},$$
as follows
\begin{equation} \label{app:eq:hub_rates_fp}
\begin{split}
\hat q&=\alpha\hat\sigma_0^2\partial_v\Upsilon(v)+\alpha(\Upsilon(v)-v\partial_v\Upsilon(v))\frac{\varepsilon_{\rm est}}{v^2}\\&\quad+\alpha\hat\sigma_0^2\partial_v\mathbb E[\erfc\chi]+\alpha\delta^2\mathbb E\left[\sigma^2\erfc\chi\right],\\
\hat m&=\alpha\frac{\Upsilon(v)}{v},    
\end{split}\qquad 
\begin{split}
	q &=(\beta_\star^2\hat m^2+\hat q)v^2\\
	m &= v\beta_\star^2\hat m,
\end{split}
\end{equation}
alongside with the equation for $v$,
\[1-\lambda v=\alpha \Upsilon(v).\]
The estimation error can be written as
\begin{equation}\varepsilon_{\rm est}=v\left(\hat\sigma_0^2+\frac{(\beta^2_\star\lambda-\hat\sigma_0^2)\lambda+\alpha\hat\sigma_0^2\partial_v\mathbb E[\erfc\chi]+\alpha\delta^2\mathbb E[\sigma^2\erfc\chi]}{\alpha\partial_v\Upsilon(v)+\lambda}\right).\label{eq:app:estsh}\end{equation}
For $\delta\to+\infty$, $\Upsilon(v)\to Y(v)$ and the equations given above coincide with Eq.~\eqref{app:eq:spC1}, as it should. The equation for $v$ specifies the relative scaling between $\alpha$, $v$, and $\varepsilon_{\rm est}\sim v$. Being $0\leq \Upsilon(v)\leq 1$, the equation for $v$ implies that, in the $\alpha\to \infty$ limit, $\Upsilon(v)\to 0$, i.e., $v\to 0$ in such a way that $\Upsilon(v)=\sfrac{1}{\alpha}+o(\sfrac{1}{\alpha})$. As $\erf\chi\to 1$ for large $\sigma$, the analytic properties of $\Upsilon(v)$ for $v\to 0$ are the same as the ones of $Y(v)$ although with different coefficients. This implies the same asymptotic scaling for the estimation error obtained for the square loss in terms of the exponent $a$ of the distribution $\varrho$ of $\sigma$.

%% file: sections/appendix/numerics.tex
\subsection{Heavy-tailed noise}
\begin{figure}
    \includegraphics[width=0.45\textwidth]{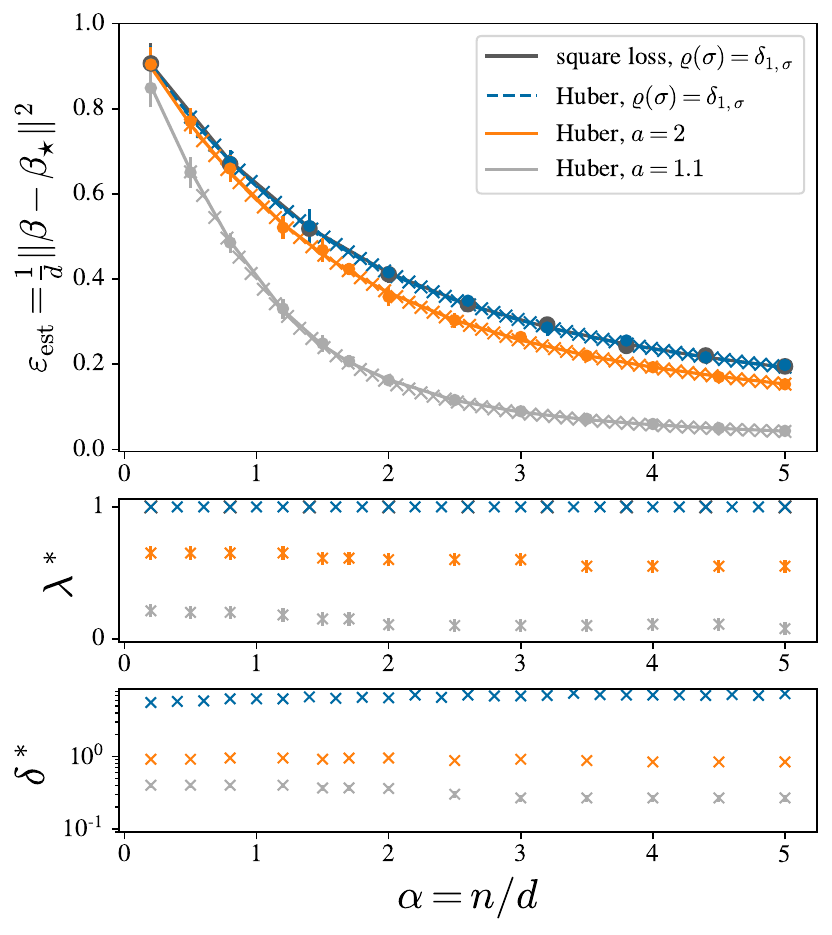}
    \quad
    \includegraphics[width=0.45\textwidth]{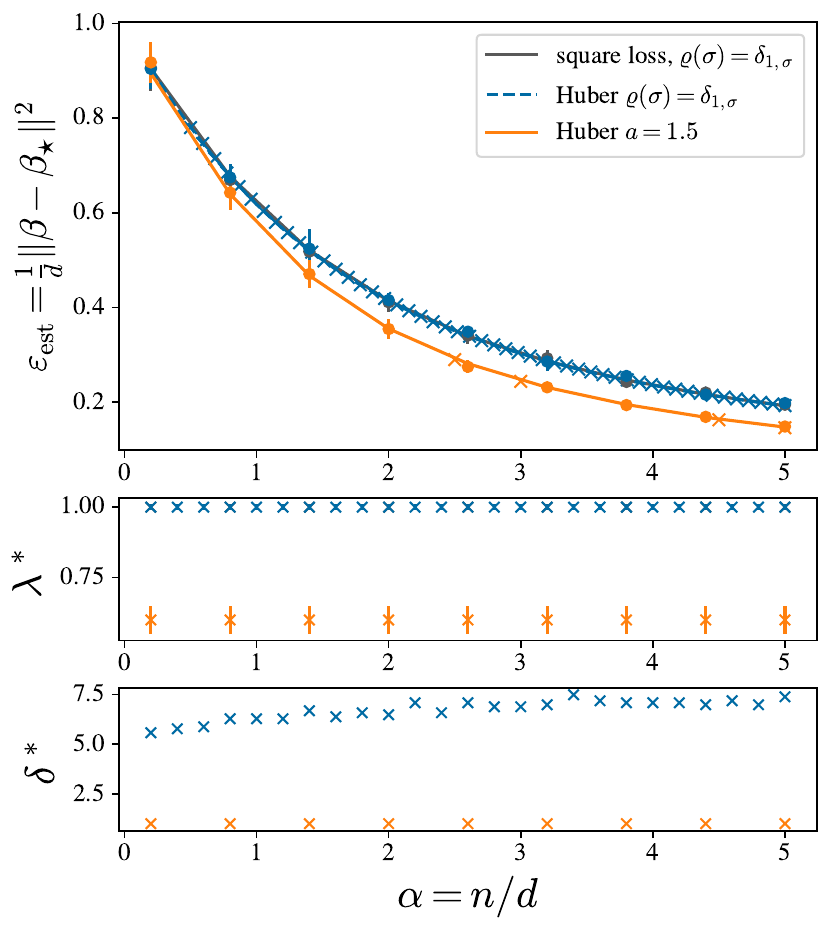}
    \caption{Gaussian covariates fully-contaminated by heavy-tailed noise with distribution $p_{\eta}(\eta)=\mathbb E_\sigma[\mathcal N(\eta;0,\sigma^2)]$, where parameters of various $\varrho(\sigma)$ are varied: inverse-gamma (\textbf{left}) and Pareto (\textbf{right}), see Table.~\ref{tab:examples}. (\textbf{Top}) Estimation error $\varepsilon_{\rm est}$ as a function of the sample complexity $\alpha = \sfrac{n}{d}$ for optimally regularised ridge regression (black), Huber with optimal location parameter and optimal regularisation (orange) and Bayes-optimal performance (crosses). Dots indicate numerical simulations averaged over 20 seeds with $d=10^3$. (\textbf{Center.}) Value of the optimal regularisation parameter $\lambda^{\star}$ for the Huber loss. (\textbf{Bottom.}) Value of the optimal location parameter $\delta^{\star}$ for the Huber loss. Both optimal values are displayed by varying the scale parameter $a$ controlling the tails of the noise distribution. }\label{sec3:fig:square_eg_rates_data_model_2_pareto_student}
\end{figure}

In this Appendix, we provide additional details about the case of $\epsilon$-contamination in the labels, as in Eq.~\eqref{eq:def:hubercontnoise}, for different $\varrho_0$ generating the contaminating noise. Fig.~\ref{sec3:fig:square_eg_rates_data_model_2_pareto_student} compares the performance of various losses for different fully-contaminated ($\epsilon_{\rm n}=1$) label noise distributions. We analyse results obtained by using for $\varrho_0$ an inverse-gamma distribution with parametrisation $a=b+1>1$ (left) and a Pareto distribution (right), see Table~\ref{tab:examples} for their full form. In all cases, the chosen parametrisations are such that unit noise variance is enforced. Note also that by taking the $a\rightarrow\infty$ the Gaussian noise limit is recovered. In our experiments, we observe the same phenomenology as in Fig.~\ref{sec3:fig:label_noise_contamination} (bottom) for all these densities which generate different noise label distributions. As long as the label noise variance is kept the same, optimally tuned Huber has the best performances as quantified by Eq.~\eqref{eq:def:esterr}. In particular, as in the cases discussed in the main text, optimally regularised Huber achieves Bayes-optimal performance, represented in the plots with crosses.

\subsection{Least Absolute Deviation estimator} \label{app:sec:LAD}
\begin{figure}
    \centering
    \includegraphics[width=0.95\textwidth]{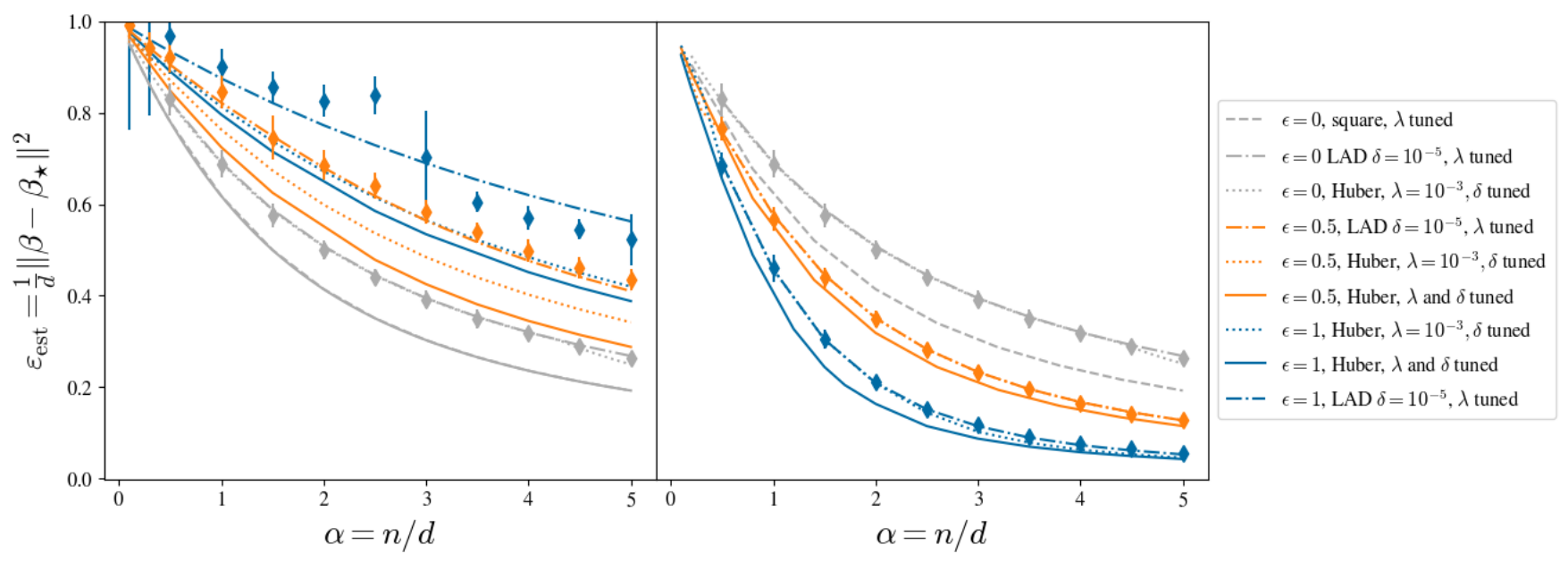}
    \caption{Label contamination including results for the LAD estimator, exactly in the same setting as Fig. \ref{sec3:fig:label_noise_contamination}. (\textbf{Left}) Infinite variance contamination. (\textbf{Right}) Finite variance contamination. Dots are experiments with dimensionality $d=10^3$, averaged over $20$ runs. }\label{app:fig:LAD}
\end{figure}

We consider the \textit{least absolute deviation} (LAD) estimator, obtained by taking the loss $\ell(y,\eta)=|y-\eta|$ in Eq.~\eqref{eq:def:erm}, with ridge regularisation. Note that this choice corresponds to a $\delta=0$ Huber loss and falls therefore within the explicit analysis above. Fig.~\ref{app:fig:LAD} compares the optimally regularised LAD estimator with the following estimators:
\begin{itemize}
    \item Optimally regularised ridge regression.
    \item Unregularised Huber with optimally chosen location parameter $\delta$.
    \item Optimally regularised Huber with optimally chosen location parameter $\delta$
\end{itemize}
We consider three label contamination cases $\epsilon_{\rm n} \in\{0,0.5, 1\}$, with both infinite and finite label variance, and Gaussian covariates, reproducing the setting analysed in Fig. \ref{sec3:fig:label_noise_contamination}. As it can be seen in Fig. \ref{app:fig:LAD}, the LAD estimator performs similarly or worse than the other estimators in all considered cases. 

\subsection{A real-world data experiment} \label{app:C3:real_data}

\begin{figure}
    \centering
    \includegraphics[width=0.9\textwidth]{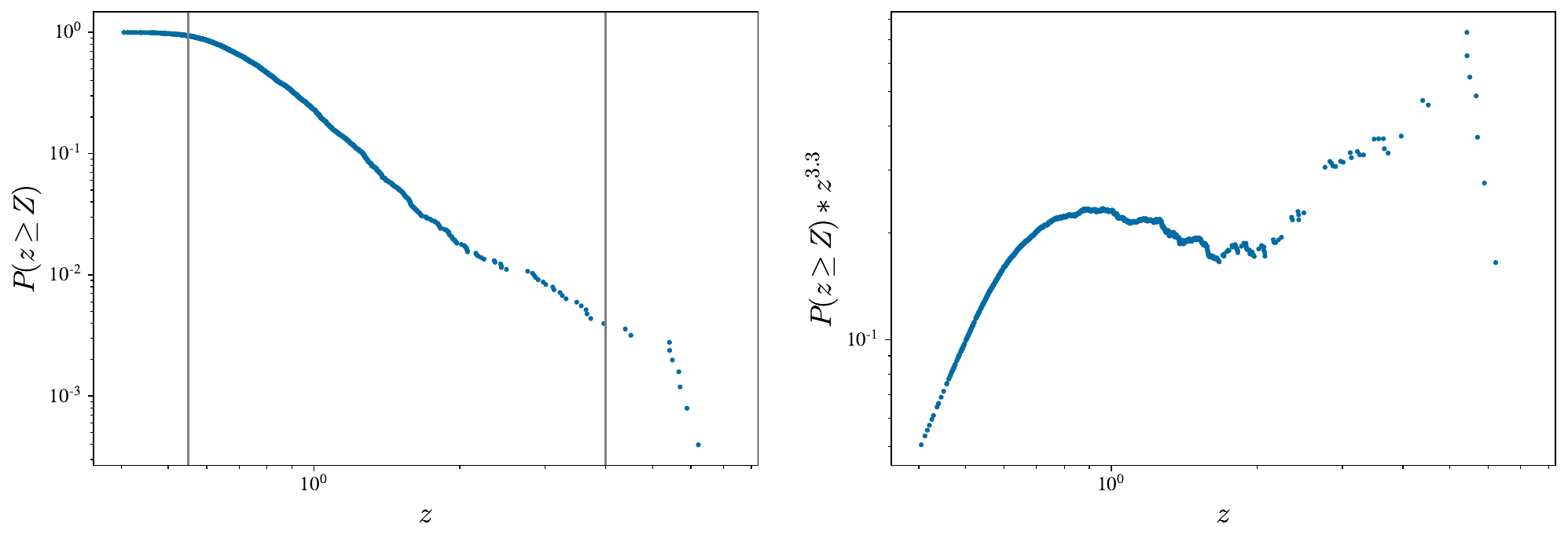}
    \caption{(\textbf{Left}) Plot of the cCDF of the norms of rotated covariates $\|\bz_i\|_2^2$, clearly exhibiting a power-law tail. The gray lines show the interval $\|\bz_i\|_2^2\in[0.55,4]$ in which the data was taken to compute the decay exponent using the \texttt{powerlaw} package \cite{powerlaw}. (\textbf{Right}) The cCDF of the norms of rotated covariates $\|\bz_i\|_2^2$, multiplied by $x^{3.3}$ to establish it as a suitable exponent that fits the power-law tail.}\label{app:fig:index_estimate}
\end{figure}

In this subsection, we present the estimation of the power-law decay index of the real data used in the experiment presented in Section~\ref{sec:covariatecont}. Such exponent is used as a parameter in the theoretical prediction of the estimation error under the assumption of heavy-tailed-distributed covariates. As anticipated in the main text, the dataset consists of a collection $\{\bx_i\}_{i=1}^n$ of $n=2519$ daily stock returns of $d=858$ stocks. Denoting by $\hat\bmu$ the empirical mean and by $\{\be_\nu\}_{\nu\in[d]}$ the basis of the empirical covariance $\hat\bSigma$ of the dataset, we introduce for each $\bx_i$ the vector $\hat\bz_i=(\hat z_i^\nu)_{\nu\in[d]}$, whose components are $\hat z_i^\nu=\sigma_0^{-1}(\be_\nu^\intercal\bx_i-\hat\mu_\nu)$, where $\sigma_0^2=\tr\hat\bSigma$, so that the empirical covariance of the new set of vectors $\{\bz_i\}_{i=1}^d$ has \textit{mean} eigenvalue $\sfrac{1}{d}$. Our primitive modeling consisted therefore in adopting as a model for the dataset a distribution in the form
\[p(\bz)=\mathbb E[\mathcal N(\bz;\mathbf 0,\sfrac{\sigma^2}{d}\bI_d)], \quad \text{where} \quad \varrho(\sigma)\coloneqq\frac{2(a-1)^{a}\exp(-\frac{a-1}{\sigma^2})}{\Gamma(a)\sigma^{2a+1}},\]
which has covariance $\bSigma=\sfrac{1}{d}\bI_d$ and whose power-law tail is determined by $a$, being $p(\bz)\sim\|\bz\|^{-2a-1}$ for $\|\bz\|\gg 1$. We have therefore one parameter only to fix, namely $a$, and we did so by estimating the power-law exponent $a$ in our dataset. Fig.~\ref{app:fig:index_estimate} shows the plot of the complementary cumulative distribution function (cCDF) $\mathbb P[\|\bz\|\geq z]\sim z^{-2a}$ of the norms of the rotated covariates, clearly showing a power-law decay in the tail of the cCDF. Moreover, we used the fitting function in the \texttt{powerlaw} package \cite{powerlaw} on norms within the interval $[0.55,4]$ whose bounds are represented in grey lines, and the decay index of the cCDF is approximated to be $2a\simeq 3.32$
. On the right, we plot cCDF multiplied by $z^{3.3}$, and the small resulting values in the bulk of the tail establish that the tail can be roughly approximated by a power-law with $2a=3.3$. 

